\documentclass[preprint,authoryear,12pt]{elsarticle}
\RequirePackage[OT1]{fontenc}

\usepackage[applemac]{inputenc}
\usepackage[french]{babel}
\usepackage{amssymb}

\usepackage{latexsym,enumerate}
\usepackage{amsmath,amsthm,amsopn,amstext,amscd,amsfonts,amssymb,euler, pdfsync,amsbsy}
\usepackage{mathrsfs}
\usepackage{fullpage}
\usepackage{graphicx}
\usepackage[usenames]{color}

\newtheorem{proposition}{Proposition}
\newtheorem{lemma}{Lemma}
\newtheorem{theorem}{Theorem}

\let\epsilon=\varepsilon

\let\phi=\varphi

\let\tilde=\widetilde

\newcommand{\field}[1]{\mathbb{#1}}
\newcommand{\R}{\field{R}}

\newcommand{\G}{{\mathcal G}}

\newcommand{\B}{{\mathcal B}}

\newcommand{\cI}{\mathcal{I}}

\newcommand{\cB}{{\mathcal B}}
\newcommand{\cC}{{\mathcal C}}

\newcommand{\beqn}{\begin{equation}}
\newcommand{\eeqn}{\end{equation}}
\newcommand\eref[1]{(\ref{#1})}

\def\I{{\field I}}

\def\lam{\lambda}

\journal{JSPI}
\begin{document}

\begin{frontmatter}
\title{Grouping Strategies and Thresholding\\
for High Dimensional Linear Models}

\author{M. Mougeot}
\author{D. Picard}
\author{K. Tribouley}

\address{ Universit\'e Paris-Diderot, CNRS LPMA, 175 rue du Chevaleret, 75013
Paris, France.}

\begin{abstract}\hspace{0.4cm}
 The estimation problem in a high regression model with
structured sparsity is investigated. An algorithm using a two steps
block thresholding procedure called GR-LOL is provided. Convergence
rates are produced: they depend on simple coherence-type indices of
the Gram matrix -easily checkable on the data- as well as sparsity
assumptions of the model parameters measured by a combination of
$l_1$ within-blocks with $l_q,q<1$ between-blocks norms. The
simplicity of the coherence indicator suggests ways to optimize the
rates of convergence when the group structure is not naturally given
by the problem and is unknown. In such a case, an auto-driven
procedure is provided to determine the regressors groups (number and
contents). An intensive practical study compares our grouping
methods with the standard LOL algorithm. We prove that the grouping
rarely deteriorates the results but can improve them very
significantly. GR-LOL is also compared with group-Lasso procedures
and exhibits a very encouraging behavior. The results are quite
impressive, especially when GR-LOL algorithm is combined with a
grouping pre-processing.
\end{abstract}

\vspace{1cm}



\end{frontmatter}

\baselineskip=18 pt

\section{Introduction}
In this paper, the following linear model is considered
\begin{align}
 Y_i&=X_i\beta+W_i,\quad i=1,\ldots,n\label{model}
\end{align}
 with a particular focus on
cases where the number $k$ of regressors $X=(X_{\bullet 1},\ldots,X_{\bullet k})$ is large compared to the
number $n$ of observations (although there is no such
restrictions).
$Y$ (respectively $W$) is denoting the $n$ dimensional observation
(respectively the error term).

 We are interested by the estimation of the parameter $\beta$ and we consider the situation where the expectation of the observation
 can be approximated by a sparse linear combination of the available regressors.
A natural method for sparse learning is $\ell_0$ regularization.
Since this optimization problem is generally NP-hard, approximate solutions
are generally proposed in practice.
 Standard approaches are
$\ell_1$ regularization, such as Lasso (see by instance
 \cite{lasso1}, \cite{Lasso3} and \cite{MeinshausenYu2009}) and Dantzig
(see \cite{dantzig}). Another commonly used approach is greedy
algorithms, such as the orthogonal matching pursuit (OMP) (see
\cite{Tropp-Guilbert}) or the iterative thresholding algorithms (see
\cite{KerkyacharianMougeotPicardTribouley}). In many practical
applications, one often knows a structure on the coefficient vector
$\beta$ in addition to sparsity. For example, in group sparsity,
variables belonging to the same group may be assumed to be zero or
nonzero simultaneously. The idea of using group sparsity has been
largely explored. For example, group sparsity has been considered
for simultaneous sparse approximation (see \cite{Wipf-Rao}) and
multi-task compressive sensing (see \cite{Ji-Dunson}) to  the tree
sparsity (see \cite{He-Carin}). Numerous applications of these types
of regularization scheme arise in the context of multi-task learning
and multiple kernel learning (see \cite{Bach}, \cite{Jenat-Bach}).
To combine sparsity with grouping, Lasso has been extended to the
group Lasso in the statistical literature by \cite{Yuan-Lin}.
 Various combinations of norms allowing grouping have been introduced as
 in \cite{Zhao-Yu}. \cite{MeierGeerBuhlmann} study the logistic
 regression model while \cite{JacobObozinskiVert} is concerning by the graph
lasso. These grouping strategies have been shown to improve the
prediction performance and/or interpretability of the learned models
when the block structure is relevant (see \cite{Koltchinski-Yuan},
\cite{Huang-Zhang}, \cite{Lounici-Sara},
\cite{ChiquetGrandvaletCharbonnier}). In
\cite{FriedmanHastieTibshirani}, Lasso and Group Lasso are combined
in order to select groups and predictors within a group.

In the sequel, we address the following program:

 $\bullet$ {\bf GR-LOL algorithm.} We investigate the theoretical
  performances of a blockwise two step thresholding algorithm. As
  LOL (standard two steps thresholding algorithm
  \cite{MougeotPicardTribouley:2012}) is a counterpart of Lasso or
  Dantzig algorithms for ordinary sparsity, we introduce here GR-LOL,
  based on the same precepts, combining an a priori knowledge of grouping.
   We establish the rates of convergence of this new procedure when
    the parameter $\beta$ belongs to a set of structured sparsity:
    the sparsity is measured by combination of $\ell_q$-between blocks
     with $\ell_1$-within blocks norms (see \eref{l1qbis}).
Although structured sparsity with overlapping groups of variables
 constitutes an important source of practical examples (hierarchical
 structure for instance), we focus in this paper on the non-overlapping case.
To emphasize the practical interest of the GR-LOL algorithm, we also
explicitly show cases where non grouping induces an accuracy loss compared to grouping.

$\bullet$ {\bf Grouping strategy.} As explained in the examples
above, in some application cases, the grouping of the predictors
occurs quite naturally or is driven by some precise requirements:
hierarchical structures or multiple kernel learning... However, in
various cases (for instance in genomic), there is no obvious
grouping at hand. In such a setting, we provide grouping strategies
which, combined with a GR-LOL algorithm, aim at improving the rates
of convergence. These grouping strategies are issued from the
following observations. Concerning the standard case (no grouping),
although  the two steps thresholding algorithms show  quite
comparable performances with Lasso and Dantzig procedures with much
less computation cost, they require theoretically more stringent
conditions on the matrix of predictors $X$, namely coherence
conditions instead of RIP-type conditions (see
\cite{MougeotPicardTribouley:2012}). In the case of structured
sparsity, this becomes surprisingly favorable, since the required
conditions -which are adaptations to the structured case of the
coherence conditions- become much more readable, and especially open
opportunities to improvements with grouping strategies. We are able
to isolate simple quantities measured on the predictors $X$ yielding
optimizing strategies to select a structure on the predictors.

$\bullet$ {\bf Practical study}. An intensive calculation program is performed to show the advantages
 and  limitations of GR-LOL procedure in several practical
aspects as well as its combination with different grouping
strategies. Based on simulations, the benefices of grouping the
predictors is compared to the non grouping case for prediction
sparse learning. We show that the way of grouping the regressors
may be critical especially when there exists some dependency
between the regressors. Using simulated data, we observe that
smart strategies of grouping strongly improve the predicted
performances.

The paper is organized as follows. In Section~\ref{section_model},
  notations and general assumptions are presented. Examples of grouping are enlightened. In
Section~\ref{GR_LOL}, the procedure GR-LOL is detailed. In
Section~\ref{results}, we state the theoretical results concerning
the performances of GR-LOL. In Section \ref{boost-section}, we first
detail explicit examples where grouping does improve the
performances, we then discuss strategies to 'boost' the rates of
convergence. The practical performances of GR-LOL  are investigated
in Section~\ref{simu} and the proofs are detailed in
Section~\ref{sectionproof}.

\section{Assumptions on the model and examples}\label{section_model}
We first introduce some notation for the predictors grouping. Next,
we state the assumptions on the model: conditions on the noise, on
the unknown parameters to be estimated and on the predictors. We end
this section with examples of models where specific grouping are
proposed.

 In the sequel, for any subset
${\mathcal I}$ of $\{1,\ldots,k\}$, $X_{\mathcal I}$ denotes the
matrix of size $n\times\#(\cI)$ composed of the columns of $X$ whose
indices are
 in ${\mathcal I}$. In the same way, $u_\cI$ is the
restriction of the vector $u$ of $\R^k$ to the vector (of dimension $\#(\cI)$) of  its coordinates with indices belonging to $\cI$. Moreover,  $\|u\|_{\cI,1}$ and
$\|u\|_{\cI,2}$  denote respectively the $l^1$-norm  and $l^2$-norm of the
restriction $u_\cI$ of $u\in\R^k$:
$$
\|u\|_{\cI,1}=\sum_{\ell\in\cI}|{u}_{\ell}|\quad\mbox{ and }\quad
\|u\|_{\cI,2}^2=\sum_{\ell\in\cI}|{u}_{\ell}|^2.
$$

\subsection{Grouping}
We consider the model \eref{model}. We consider a partition
$\G_1,\ldots,\G_p$ of the set $\{1,\ldots,k\}$ of the indices of the
regressors. For any $j$ in $\{1,\ldots,p\}$,
 $t_j=\#(\G_j)$ denotes the cardinal of the group $\G_j$.
We decide to subdivide the $k$ predictors into $p$ ($p\leq k$)
groups of variables $X_{\G_1},\ldots,X_{\G_p}$, according to this partition.
  Following this subdivision, for each $\ell=1,\ldots,k$,
 the predictor $X_\ell$ is now registered as $X_{(j,t)}$ where
\begin{itemize}
\item $j\in\{1,\ldots,p\}$ is the index of the group $\G_j$ where
the index $\ell$  belongs,
\item  $t=r_j(\ell)\in\{1,\ldots,t_j\}$ is the rank of $\ell$ inside the group
$\G_j$.
\end{itemize}
The notation $\ell=(j,t)$ is used all along the paper. The group of
indices $\G_j$ is then identified with $\{(j,t)\mbox{ for
}t=1,\ldots,t_j\}$. The index $t$ will sometimes in the sequel be assimilated to a 'task index' in analogy to the forthcoming example \ref{multitask-section}.

\subsection{Assumptions}
\subsubsection{Homogeneousness condition for the predictors }
To take into account the natural inhomogeneity of the data, we
define a normalizing constant $n_\ell$ depending on $\ell\in
\{1,\ldots,k\}$. It appears naturally as a 'normalizing constant'
through the forthcoming assumption \eref{norm}. Setting $\tilde
X_{i\ell}=X_{i\ell}/\sqrt{n_\ell}$ for any observation
$i=1,\ldots,n$, the model becomes
\begin{align}
 Y&=\tilde X\;\alpha+W\label{modeltilde}
 \end{align}
where
$$\alpha_\ell= \sqrt{n_\ell}\beta_\ell\quad \mbox{ for any }\; \ell=1,\ldots,k$$
In the sequel we assume that there exists a sequence $v_n$  and
constants $0<a<b$ such that for any $\ell \in \{1,\ldots,k\}$, we
get
\begin{equation}
 (A1):\quad a \,v_n\le n_\ell\le b\, v_n.\label{homo}
 \end{equation}
The quantity $\nu_n$ is important because it drives the rates of
convergence of our algorithm.
\subsubsection{Conditions on the predictors}\label{condpredic}
Denote $\Gamma=\tilde X^t\tilde X$ the Gram matrix of $\tilde X$ and
$\Gamma_\cI=\tilde X_\cI^t\tilde X_\cI$ the Gram matrix of $\tilde
X_\cI$. Observe that $\Gamma_{\ell\ell^\prime}$ is the scalar
product between two predictors $\widetilde{X_\ell}$ and
$\widetilde{X_{\ell^\prime}}$ and define the coherence of the Gram
matrix
\begin{align}
\gamma:=\sup_{(\ell,\ell^\prime)\in\{1,\ldots,k\}^2,\ell\not =
\ell^\prime} \left|\Gamma_{\ell\ell^\prime}\right|\label{cohe}
\end{align}
 Recall that each $\ell$ in $\{1,\ldots,k\}$ is registered as
a pair of indices $(j,t)$ where $j$ is the index of the group and
$t$ the rank of $\ell$ inside the group and then
$$
\ell\not =\ell^\prime\Longrightarrow \quad\left\{\begin{array}{ll}
t\not =t^\prime &\ell\mbox{ and }\ell^\prime\mbox{ are not observed at the same 'task'}\\
\mbox{ or }&\\
 t=t^\prime, j\not =j^\prime\;&\ell\mbox{ and }\ell^\prime\mbox{ are observed at the same 'task' but in different groups}\end{array}\right.
$$
We split $\gamma$ as $\gamma_{BT}\vee\gamma_{BG}$ where
\begin{align}
\gamma_{BT}:=\sup_{(j,j^\prime)\in\{1,\ldots,p\}^2}\sup_{t\in\{1,\ldots,t_j\},t^\prime\in\{1,\ldots,t_{j^\prime}\},
t\not =t^\prime}
\left|\Gamma_{(j,t)(j^\prime,t^\prime)}\right|\label{btask}
\end{align}
and \begin{align}
\gamma_{BG}:=\sup_{(j,j^\prime)\in\{1,\ldots,p\}^2,j\not =j^\prime }
\sup_{t\in\{1,\ldots,t_j\wedge t_{j^\prime}\}}
\left|\Gamma_{(j,t)(j^\prime,t)}\right|.\label{bgroup}
\end{align}
For any
 subset $\cI$ of the set of indices $\{1,\ldots,k\}$, let $\tau(\cI)$ and $r(\cI)$
   be the following indicators
\begin{align}\label{tauI}
\tau(\cI):=\#(\cI)\;\gamma_{BT}+\#(\{j,\; \exists t,\; (j,t)\in\cI\})\;\gamma_{BG}.
\end{align}
\begin{align}\label{tauI}
r(\cI):=\#(\cI)\;\gamma_{BT}^2+\#(\{j,\; \exists t,\;
(j,t)\in\cI\})\;\gamma_{BG}^2.
\end{align}
In particular,  for any $j\in\{1\ldots, p\}$, we define
\begin{align*}
\tau_j=\tau(\G_j)=t_j\;\gamma_{BT}+\gamma_{BG}\quad\mbox{ and }\quad
r_j=r(\G_j)=t_j\;\gamma_{BT}^2+\gamma_{BG}^2
\end{align*}
as well as
\begin{align}\label{tau*}
\tau^*=\max_{j=1,\ldots,p}\tau_j=t^*\;\gamma_{BT}+\gamma_{BG}
\quad\mbox{ and }\quad
r^*=\max_{j=1,\ldots,p}r_j=t^*\;\gamma_{BT}^2+\gamma_{BG}^2.
\end{align}
where $t^*=\max_{j=1,\ldots,p}(t_j)$.

Let us state now the assumptions on the regressors $X$. First, we
assume that the columns of the matrix $\tilde X$ are normalized:
\begin{align}
(A2):\Gamma_{\ell\ell}=1\quad \mbox{ for any
}\,\ell=1,\ldots,k.\label{norm}
\end{align}
 Second, we assume that
\begin{equation}
(A2'):\quad \tau^* \le \nu.\label{mini}
\end{equation}
for some $\nu$ given in $]0,1[$. Observe that under $(A2)$, we
obviously have $r^*\leq \tau^*$.
\subsubsection{Conditions on the unknown regression parameters}
Assume that there exist $q\le 1$ and $M,M^\prime>0$ such that
\begin{equation}\label{l1qbis}
(A3):\quad\sum_{j=1}^p\|\beta\|_{\G_j,1}^{q}\le (M^\prime)^{q} \quad\mbox{  or
equivalently } \quad\sum_{j=1}^p\|\alpha\|_{\G_j,1}^{q}\le
M^{q}\;v_n^{q/2}.
\end{equation}
\subsubsection{Conditions on the noise}
Finally, we assume
\begin{equation*}
(A4):\quad  W \mbox{ is a vector of i.i.d. variables } {\mathcal
N}(0,\sigma^2).
\end{equation*}
 Notice that the Gaussian distribution assumption may
be replaced without modifications by a sub-Gaussian distribution
with zero mean and variance $\sigma^2$.
\subsection{Specific models. Examples.}
\subsubsection{No-group case}

One specific case of our modeling is when $\G_j=\{\,j\,\}$ for any
$j\in\{1,\ldots,k\}$: the no-group setting which corresponds to
$p=k$. Here, the predictors are generally normalized by the number of
observations $n_\ell=n$ and the homogeneousness condition
(\ref{homo}) is ordinary satisfied for $v_n=n$. Moreover, $\gamma_{BT}=0$ and
$$\gamma_{BG}=\sup_{(\ell,\ell^\prime)
\in\{1,\ldots,k\}^2,\ell\not=\ell^\prime}|\Gamma_{\ell\ell^\prime}|
$$
is the coherence of the matrix $\Gamma$. We get $\tau^*=\gamma_{BG}$
and Condition (\ref{mini}) becomes $\gamma_{BG}\leq\nu$. Note that a similar
 condition is used in \cite{KerkyacharianMougeotPicardTribouley}
or \cite{MougeotPicardTribouley:2012}. The regularity conditions in
this case sum up to a $l_q$ condition on the parameter vector
$\beta$.

\subsubsection{Multi-task case}\label{multitask-section}

An interesting case where many conditions find direct interpretation is  the multi-task
regression model defined by the pile of $T$ linear models:
\begin{eqnarray}\label{multitask}
\left\{ \begin{array}{ll}
Y_1&=X_1\beta_1+W_1\\
Y_2&=X_2\beta_2+W_2\\
&\ldots\\
Y_T&=X_T\beta_T+W_T
\end{array}\right.
\end{eqnarray}
Here $X_1$, ..., $X_T$ are  $n_0\times p$ design matrices
and  $W_1,\ldots,W_T$ are (independent) error terms. This modeling is used (for
instance) to introduce a time variation: the target variable $Y$ and
the predictors $X_1,\ldots,X_p$ are observed on $T$ different periods of time.
We prefer the term task to time not to induce confusion with
the 'observation times' $i$. For each task the observation consists in a vector
 $Y_t$ of size $n_0$, analyzed on the matrix of predictor $X_t$.
Model (\ref{multitask}) can be globally rewritten as Model \eref{model} with
$$n=n_0T\quad\mbox{ and }\quad k=pT,$$
 the design matrix $X$ being block diagonal with blocks $X_1,\ldots,X_T $
and   $$\beta=(\beta_1^t,\ldots,\beta_T^t)^t\quad\mbox{ and }\quad
Y=(Y_1,\ldots,Y_T)^t.$$ We obviously have $n_\ell=n_0$ for any
$\ell=1,\ldots,k$ and the normalization condition (\ref{homo}) is ordinary
satisfied for $v_n=n_0$.
 Notice that
 the different groups of indices, $\G_j=\{(j,t),\; t=1,\ldots ,T \}$
  for $j=1,\ldots,p$ have all the same size
 $T$; the index $j$ points out the predictor $X_j$ for $j=1,\ldots,p$
 and the index $t$ is an indicator of the task of observation for
 $t=1,\ldots,T$. Thanks to the block structure of the matrix $X$,
 the predictors are obviously orthogonal as soon as the
  tasks are different; even the same variables observed
 at different tasks are orthogonal. We deduce that
$\gamma_{BT}=0$. Moreover, denoting $\Gamma_1,\ldots,\Gamma_T$ the
sequence of Gram matrices associated to the $T$ models given in
(\ref{multitask})
$$\gamma_{BG}=\max\left(\sup_{j\not=j^\prime}|(\Gamma_1)_{jj^\prime}|,\,
\ldots \,,\sup_{j\not=j^\prime}|(\Gamma_T)_{jj^\prime}|\right)
$$
and Condition (\ref{mini}) becomes $\gamma_{BG}\leq\nu$.

This example is especially emblematic. In this context, the rank $t$
in the group  $\G_j$ is easily interpretable as a task. As well,
condition (A3) is quite  realistic since the coefficients $\beta_{(j,t)}$ on
 the predictor $X_{(j,t)}$ can be assumed to slowly vary with the task.
 Furthermore, the separation introduced in subsection \ref{condpredic}
 between $\gamma_{BT}$ and $\gamma_{BG}$, which, in an implicit way assumes
  in condition (A2') that $\gamma_{BT}$ is a smaller quantity, naturally finds
  its interpretation here (since it is 0).

\section{GR-LOL: Grouping Research for Leaders}\label{GR_LOL}

Let us now describe the steps of our procedure. Once for all, we fix
the constant $\nu$ which is a quantity linked to the precision of
the procedure; take for instance $\nu=1/2$.

\begin{description}
\item[Compute a bound for the number of leaders.]
Form $\Gamma=\tilde X^t\tilde X$ and compute
$\gamma_{BT},\gamma_{BG}$ as defined in (\ref{btask}) and
(\ref{bgroup}). Deduce $\tau^*=t^*\,\gamma_{BT}+\gamma_{BG}$ and
$N^*=\nu\, (\tau^*)^{-1}$  (see Definition (\ref{tau*})) .

\item[Search for the leaders.]
Form $ R_\ell=\sum_{i=1}^nY_i\;\widetilde{X_{i\ell}}$ (
rewritten as $R_{(j,t)}$ to take into account the group number) and compute, for any group
$\G_j,\,j=1,\ldots,p,$ the quantity
 $$\rho_j^2=\sum_{t=1,\ldots,t_j} R_{(j,t)}^2:=
\| R\|_{\G_j,2}^2\,.
$$
$\rho_j^2$ is an indicator of performance for the predictors whose  indices are in the
group $\G_j$  to explain the target variable $Y$. Next, we consider
the groups for which this indicator is high. More precisely, the
sequence $\rho_j^2$ is sorted:
$$
\rho^2_{(1)}\geq \ldots \geq \rho^2_{(j)}\geq \ldots \geq
\rho^2_{(p)}
$$
and the group-leaders are the groups of predictors with group-indices in
 $j\in\cB$ where
\begin{eqnarray}\label{leader}
\cB=\,\left\{ j=1,\ldots,p,\; \rho_{(j)}^2\geq
\left(\rho_{(N^*)}^2\vee \lam_n(1)^2\right)\right\}
\end{eqnarray}
where $\lam_n(1)$ is a first tuning parameter. Denote
$\G_\cB=\cup_{j\in\cB}\G_j$. Notice here that in the case where
$\lam_n(1)^2>\rho^2_{(1)}$, the leader indices set $\cB$ is empty
and our final estimate for $\beta$ is zero.

Observe also that $\#(\cB)\leq N^*$ and $\#(\G_\cB)\leq t^*N^*$ implying that
$$
\tau(\G_\cB)\leq N^*(t^*\,\gamma_{BT}+\gamma_{BG})=\nu.$$

\item[Regress on the leaders.]
We now  perform the
OLS on the block-leaders:
$$\hat\beta(\cB)=\hbox{Argmin}_{u} \|Y-X_{\G_\cB}u\|^2=
[X_{\G_\cB}^tX_{\G_\cB}]^{-1}X_{\G_\cB}Y.$$

We then obtain  the preliminary estimate $\hat\beta$ defined by
$$
 \hat\beta_{\G_\cB}=
\hat\beta(\cB)\quad \mbox{ and }\quad \hat\beta_{\G_\cB^c}=0$$

\item[Block thresholding]

We apply the second thresholding on the resulting estimated coefficients:
$$\forall \ell=(j,t)\in\{1,\ldots,k\},\quad\hat\beta^*_\ell=\hat\beta_\ell\;
\I\{\;\|\hat{\beta}\|_{\G_j,2}\ge
\frac{\lam_n(2)}{\sqrt{n_l}}\,\}\,$$ where $\lam_n(2)$ is the second
tuning parameter.
\end{description}

\section{Results}\label{results}
In this section, we provide a result on the convergence rate of
GR-LOL procedure for a quadratic error on the estimation on the
$\beta$ coefficients on the regression model when the input
parameters  $\lam_n(1),\lam_n(2)$ are properly chosen.

\subsection{Rates of convergence}
The proof of the following theorem is given in Section
\ref{sectionproof}.
\begin{theorem}\label{JAIME} Fix $\nu$ in $(0,1)$ and assume that
$A(1),A(2),A(2'),A(3)$ and $A(4)$ are satisfied. Put
\begin{align}\label{lam0}
\lambda^*=\sigma\left( M^2\,v_n\,r^*\,\vee \,(t^*\vee \log
p)(1+\nu)\right)^{1/2}.
\end{align}
Choose the thresholding levels $\lam_n(1),\lam_n(2)$ such that
$$\lam_n(2)=c_2\,\lam^*\quad\mbox{ and }\quad \lam_n(1)=c_1\,\lam^*\vee
 \left(2\,M^{ }\,v_n^{1/2}\,\tau^*/\nu\right)$$
 for $$c_1>c_2,\;c_2\geq 5\sqrt{\kappa},\;
c_1\geq (4+\nu^{-1/q})$$ and
$$\kappa=(1-\nu)^{-1}\vee 4(1-\nu)^{-2}\vee
2(2\nu^2-\nu+3)(1-\nu)^{-3}.$$ There exists a
 positive constant $C$ (depending on $c_1,c_2,\nu$ and $M$) such that
\begin{align*}
E\|\hat\beta^*-\beta\|_2^2&\leq C\,\left[\frac{
t^*\vee\log{p}}{v_n}\vee (\tau^*)^2\right]^{1-q/2}
\end{align*}
as soon as
\begin{align*}
p&\leq c_a^{-1}\,v_n^{q/2}\,(\lam^*)^{-q}\;
\exp\left(c_b\,(\lam^*)^2(1\wedge(r^*)^{-1}) \right)
\end{align*}
where
$$c_a=M^{-q}\left(\frac 94 c_1^2\vee 3\kappa\right)\;\mbox{ and }\;
c_b=\frac{c_1^2}{64(1+\nu)}\wedge \frac{c_2^2}{192\kappa}.
$$
\end{theorem}
\subsection{Comments}
It worthwhile to notice that Theorem \ref{JAIME} rather clearly
identifies the key features needed for our procedure to be sharp.
Basically, it is depending on the structured sparsity as well as the
size of the groups and the correlation structure within task and
groups.
\begin{description}
\item[Structured sparsity] Concerning the structured sparsity of the
coefficients, condition \eref{l1qbis} reflects overall an
homogeneousness inside the groups as well as a small number of
'significant' groups. As is illustrated in Section \ref{gr/non-gr},
the algorithm has better rates if the large coefficients are
gathered in the same groups, instead of being  scattered in
different groups.
\item[Size and correlation inside the groups] A key quantity is
$\tau^*=t^*\gamma_{BT}+\gamma_{BG}$. In particular, this quantity
gives clear some indication to optimize the procedure when the
structure is not a priori given by the problem. This is detailed in
the following section.
\end{description}

\subsection{A specific example: No-group case}
 In the no-group case,  the performances stated in the previous
theorem are similar to those achieved by the standard LOL procedure
studied in \cite{MougeotPicardTribouley:2012}. Actually, in the
no-group case, $p=k$;  recall that $v_n=n$, $ n_\ell=n$, $t^*=1$,
$\gamma_{BT}=0$ and that
  $\tau^*=\gamma_{BG}$ is the coherence of the
matrix $\Gamma$. Observe also that in this case $r^*=(\tau^*)^2$.
Condition $A(3)$ (see (\ref{l1qbis})) is here the usual $l_q$
condition. Applying Theorem \ref{JAIME}, we choose
$$\lam_n(1)=c_1\left(\,n^{1/2}\gamma_{BG}\vee \sqrt{\log{k}}\right)\quad \mbox{ and }\quad
\lam_n(2)=c_2\left(\,n^{1/2}\gamma_{BG}\vee \sqrt{\log{k}}\right)$$
for constants $0<c_2<c_1$ large enough and we get
\begin{align*}
E\|\hat\beta^*-\beta\|_2^2&\leq C\left(\gamma_{BG}^2\vee\frac{
\log{k}}{n}\right)^{1-q/2}
\end{align*}
under the condition
\begin{align*}
k&\leq c\left(n^{q/2}\,(n\gamma_{BG}^2\vee\log{k})^{-q/2}\right)\;
\exp\left(c\,n\gamma_{BG}^2\vee\log{k}\right)
\end{align*}
which writes as a lower bound for the constants above when
$n\,\gamma_{BG}^2\sim \log k$. There is no limitation on $k$ except
$\log k/n\le C$. In this case, the rate is minimax.
\subsection{An more interesting example: Multi-task case}
In the multi-task case, we observe $n_0$ observations issued from
$p$ variables on $T$ tasks units. We have
$$v_n=n_0,\quad n_\ell=n_0,\quad t^*=T,\quad \gamma_{BG}=0$$
and
$\tau^*=\gamma_{BG}$
is the maximum of the coherences associated to the different Gram sub-matrices
$\Gamma_1,\ldots,\Gamma_T$. As previously, we get $r^*=(\tau^*)^2$.
Choosing
$$\lam_n(1)=c_1\left(\,n_0^{1/2}\gamma_{BG}\vee \sqrt{T}\vee \sqrt{\log{p}}\right)\quad \mbox{ and }\quad
\lam_n(2)=c_2\left(\,n_0^{1/2}\gamma_{BG}\vee \sqrt{T}\vee
\sqrt{\log{p}}\right)$$ for constants $0<c_2<c_1$ large enough and
we get
\begin{align*}
E\|\hat\beta^*-\beta\|_2^2&\leq C\left(\gamma_{BG}^2\vee\frac{
T}{n_0}\vee\frac{ \log{p}}{n_0}\right)^{1-q/2}
\end{align*}
under the condition
\begin{align*}
p&\leq c\left(n_0^{q/2}\,(\,n_0\gamma_{BG}^2\vee T\vee
\log{p})^{-q/2}\right)\; \exp\left(c\,(\,n_0\gamma_{BG}^2\vee T\vee
\log{p})\,\right)
\end{align*}
yielding a lower bound for the constants here above when
$n_0\,\gamma_{BG}^2\sim T\sim \log p$. Observe there is no
limitation on $p$.

\subsection{'Minimaxity', comparisons}
In this part, we use Theorem \ref{JAIME} to evaluate the quality of
our procedure in various cases.
\begin{description}
\item[Minimax-no group] In the no-group case, minimax bounds are known (see
\cite{raskutti2011}) and our procedure achieves this bounds as soon
as $\gamma_{BG}=\tau^*\le O\left(\sqrt{\log k/n}\right)$.
\item[Still minimax when grouping] For any $q\le 1$, we obviously have
$$\sum_{j=1}^p\|\beta\|_{\G_j,1}^{q}\le
\sum_{j=1}^p\|\beta\|_{\G_j,q}^{q}=\sum_{\ell}|\beta_\ell|^{q}.$$
Hence, as soon as $\tau^*\leq O\left(\sqrt{\log k/n}\right)$ which
is satisfied for instance if
$$\gamma_{BT}\simeq 0,\quad\gamma_{BG}\le
O\left(\sqrt{\log k/n}\right)\mbox{ and }t^*\le \log k,$$ the GR-LOL
procedure is still minimax using again the lower bound given in
\cite{raskutti2011}.
\item[Wavelet coefficients] Let us consider the standard case of the signal
model where $k\le n$ and where the $\beta$'s are the wavelet
coefficients of the unknown signal. Observe that the condition
 $\|\beta\|_{q}\le M$ for $q\le 1$
corresponds to belonging of the signal to a ball of the Besov space
$B^{1/q-1/2}_{q,q}$.
 Hence Theorem \ref{JAIME} proves that GR-LOL is
 minimax for any grouping strategy such that
 $$\tau^*\le O\left(\sqrt{\log k/n}\right)\quad\mbox{ and }\quad t^*\le
 O\left(\log k\right).$$ This is
 an extension of the block thresholding strategies which are
 generally performed with  blocks  chosen inside each multiresolution level
 (see for instance, among many others \cite{HKP}, \cite{Cai-Zhou09}).
\item[Comparison with other structured sparsity conditions]
Our conditions involving simple correlation quantities on the
regressors are quite difficult to compare with more involved
conditions of geometric nature, as in \cite{Lounici-Sara} or of
structured sparse coding  nature as in \cite{Huang-Zhang} for
instance. Let us just mention that these conditions are  very likely
to be stronger than other ones, as it is the case in the no-group
case compared to RIP conditions. However they have the advantage of
being checkable on the data and they are readable enough to give
directions to optimize the procedure. This point is developed in the
sequel providing an algorithm to determine the groups.
\end{description}

\section{Boosting the rates using grouping }\label{boost-section}
Generally in structured sparsity frameworks, the grouping is coming
from the data, as it is the case for instance in the multitask case.
However, in many situations there is no indication for such a
'natural' grouping. Our purpose is to explain how proceed for
boosting the rates using grouping. We investigate different ideas
for grouping strategies and in subsection \ref{BRG}, a new grouping
(auto driven) procedure called "Boosting Rates Gathering" is
provided. To better introduce the BRG algorithm, we first detail an
example explaining what gain can be expected by a suitable grouping
and to what extend.

To simplify (but with obvious generalization), we assume that
$\nu_n=n$ in this section.

\subsection{Grouping versus non-grouping}\label{gr/non-gr}
 Consider a model such that the Gram
matrix $\Gamma$ is such that $\gamma\ge \sqrt{(\log{k})/n}$ (which
is the standard case).

$\bullet$ {\bf Use GR-LOL.} First, assume that the grouping is such
that $\gamma_{BT}\le \gamma/{t^*}$ and $ \gamma_{BG}\le \gamma$.
Assume in addition that $ t^*\le c[\gamma^{-1}\vee n\gamma^2]$ for
some positive constant $c<1$. We see below that these conditions can
automatically be ensured by the following BRG algorithms.

Consider the case where
$$\beta_\ell=\left\{\begin{array}{ll}
\gamma &\mbox{ if }\ell \in {\mathcal
G}_1\cup\,\ldots\,\cup{\mathcal G}_{\lfloor (\gamma
t^*)^{-q}\rfloor}\\
0&\mbox{ else }\end{array}\right.
$$
So we have $\#\{\ell,\,\beta_\ell\not =0\}\leq t^*\,\lfloor (\gamma
t^*)^{-q}\rfloor$. Since
$$\sum_{j\le \lfloor (\gamma
t^*)^{-q}\rfloor }(t^*\gamma)^q\le (\gamma
t^*)^{-q}(t^*\gamma)^q=1,$$
 Condition (A3) is then fulfilled with
$M=1$, Then applying Proposition \eref{boost}, the predicted error
is bounded by $C\,\gamma^{2-q}$.

$\bullet$ {\bf Use LOL.} Second, we use LOL (corresponding to GR-LOL in the no-group case) and we denote $\hat\beta^@$ the estimate obtained using this second algorithm. Since
$\tilde X=n^{-1/2}\,X$, recall that
\begin{align*}
R_\ell&=\sum_{i=1}^n Y_i\tilde X_{i\ell}=\sum_{i=1}^n
\left(X\beta+W\right)_i\tilde X_{i\ell} =\sum_{i=1}^n \left(
\sum_{\ell^\prime=1}^k \beta_{\ell^\prime}X_{i\ell^\prime}\tilde
X_{i\ell}\right)+\sum_{i=1}^n \tilde X_{i\ell}W_i\\
& = n^{1/2}\beta_\ell\, \Gamma_{\ell\ell}+
n^{1/2}\,\sum_{\ell^\prime=1,\ldots,k,\ell^\prime\not=\ell,\beta_{\ell^\prime}\not=
0}\beta_{\ell^\prime} \Gamma_{\ell\ell^\prime}+\sum_{i=1}^n \tilde
X_{i\ell}W_i.\end{align*} We deduce
\begin{align*}
R_\ell &= a_\ell+b_\ell +\xi_\ell
\end{align*}
where $$|a_\ell|\leq n^{1/2}\gamma,\qquad |b_\ell|\le
n^{1/2}\gamma^2\,\left(t^* \lfloor (\gamma
t^*)^{-q}\rfloor\right)\le 2n^{1/2}\gamma $$ and
 $\xi_\ell$ is distributed as a centered gaussian distribution of
variance $1$. Choose now in Theorem \ref{JAIME}, $\lam_n(1)\ge
5n^{1/2}\gamma$ (this choice is compatible with the assumptions in
there). Then we get, for any index $\ell$ associated with a non zero
coefficient $\beta_\ell$
\begin{align*}
P(|R_\ell|\le \lam_n(1))&\ge 1-P(|R_\ell-ER_\ell|> 2
n^{1/2}\gamma)\ge
1-\;\exp\left(-\frac{(2\sqrt{n}\,\gamma)^2}2\right)
\end{align*}
which can be bounded below by $0.5$ for $\sqrt{n}\,\gamma$ larger than an
absolute constant. Since
$$
|\hat\beta^@_\ell-\beta_\ell|\geq
|\hat\beta^@_\ell-\beta_\ell|\,\I_{\ell\not\in {\mathcal
B}}=|\beta_\ell|\,\I_{\ell\not\in {\mathcal B}}=\gamma
\,\I_{|R_\ell|\le \lam_n(1)},
$$
we deduce
$$E\|\hat\beta^@-\beta\|_2^2\ge
0.5\,\gamma^2 \,\left(t^*( \lfloor (\gamma t^*)^{-q}\rfloor\right)$$
and the predictor error is always larger than  $0.5\,(t^*)^{1-q}\,
\gamma^{2-q} $. So the prediction using grouping gives an average
error smaller by  a factor of $(t^*)^{1-q}$ which can rapidly be
substantially large when $t^*$ itself grows.

Observe also that the first procedure takes benefit of the fact that
the  'big' (here the non zero) $\beta$'s are 'gathered' in the same
groups. If instead, we have a configuration with the same final
number of $\beta$'s, all equal to $\gamma$, but scattered all in
different groups, then condition (A3) is no longer satisfied and the
group procedure achieves a lower rate. Actually a closer look at the
proofs shows that the rate is the same as obtained by the LOL
procedure.

\subsection{Gathering}\label{Gather}
A natural idea coming from the example above is to 'gather' in the
same group the indices $\ell$'s with $R_\ell$ substantially big or
of the same size. This obviously helps  to decrease the number of
groups which is an important issue. Natural ways to proceed are the
gathering procedures below.
\begin{itemize}
\item  ({\bf GGa}) Gathered Grouping with absolute correlation: this procedure gathers,
 in each group, variables exhibiting
similar absolute value  $|R_\ell|$ of the correlation coefficients with
the target $Y$. The $p$
different groups are then successively filled by using the ordered
indices:
$$
\G_1=\{\,(1),\ldots,(\lfloor k/p\rfloor)\,\},\quad \ldots\quad
,\G_p=\{\,(k-\lfloor k/p\rfloor),\ldots,(k)\,\}
$$
where $(\ell)$ denotes the index associated to the ranking quantity
$|R_{(\ell)}|$.
\item  ({\bf GGc}) Gathered Grouping with  correlation: it is the same
procedure as ({\bf GGa})
but using the $R_\ell$'s
instead of the absolute value $|R_\ell|$'s.
\end{itemize}
In view to explore in practice the benefit of these grouping
strategies (see the next section), we also introduce
\begin{itemize}
\item  ({\bf GGr}) Random  Gathered Grouping: this procedure gathers, in
each group, $k/p$ variables randomly chosen among the $k$ regressors.
\end{itemize}

\subsection{Taking into account the coherence and the size $t^*$}
If we look at the convergence results of Theorem \ref{JAIME} in view
to boost the rates, we observe that not only the structured sparsity
is important but also that the following quantity has to be
optimized

\begin{equation}
[\sqrt{\frac{t^*\vee\log p}n}\vee \{t^*\;\gamma_{BT}+\gamma_{BG}\}].\label{tobalance}
\end{equation}
Looking at this quantity gives some indications for choosing a
procedure. First $t^*$ has to be smaller than $\log p$ if possible.
This obviously induces to choose balanced groups. Looking now at the
quantity $ \tau^*= t^*\;\gamma_{BT}+\gamma_{BG}$ indicates that the
rates would benefit of choosing groups in such a way that
$\gamma_{BT}$ is as small as possible. As a consequence
$\gamma_{BG}$ is equal to  the maximal correlation $\gamma$. This
observation gives rise to the following strategy. Divide the columns
of $X$ into two sets: $S_1$ of the items which are highly
correlated, $S_2$ for the remaining,
 weakly correlated.
Put $S_1$ all in 'Task' number 1: we ensure then that $\gamma_{BT}$
is less than the maximal correlation within $S_2$ while
$\gamma_{BG}=\gamma_{max}=\gamma$. Another way to describe this is
that  each  columns of $S_1$ is the first point of a new group. This
induces in the sequel the name of 'delegate'.

It now remains to answer the two questions: how to choose the number
of groups (cardinal of $S_1$) and how to fill up the groups after
the choice of its delegate. The answers to these questions are
obtained by balancing the quantities in (\ref{tobalance}), and then
using the gathering principle. A final remark is that  the quantity
$\gamma$ is generally a leading term. Let us now be more precise and
describe BRG the procedure (Boosting Rates Gathering)
\subsection{BRG (Boosting Rates Gathering)\label{BRG}}
\subsubsection{Determination of the number $p^*$ of groups\label{numberofgroups}}
\begin{figure}[h]
\begin{center}
\includegraphics[width=7.7cm,height=4cm]{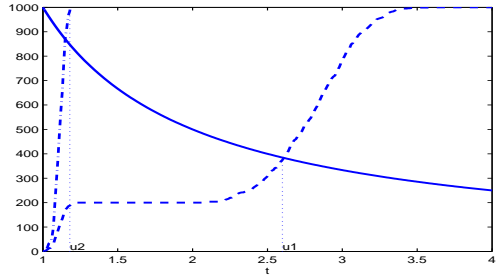}
\caption { $X$-axis: Common size $t_1$. $Y-$axis: number $p$ of
groups. Solid line: $g(u)=k/u$. Dashed line: $p(u)$ for  $\rho=0.5,\pi=20\%$ (see simulation part). Dot dashed line: $p(u)*log p(u))$.
Dot lines: corresponding $u_1$, $u_2$ positions. $n=200$, $k=1000$, $SNR=5$. We observe $u_2 < u_1$. \label{fg}}
\end{center}
\end{figure}

This is the first step of the BRG procedure. Since we choose to have
balanced groups, it is equivalent to determine the number of groups
$p$ or the average size $t^*=k/p$ of the groups. Let us consider the
following curves $y=g(u)$ and $y=p(u)$ defined for $u$ in
$[1,\infty[$,
$$g(u)=k/u\;\mbox{ and }\; p(u)=\#\left\{\ell \in\{1,\ldots,k\},\,\exists
\ell^\prime\in\{\ell+1,\ldots,k\}\mbox{ such that } |
\Gamma_{\ell\ell^\prime}|>\gamma/u \right\}.$$ These curves
intersect at a point $u_1$ as illustrated in Figure \ref{fg}.
Observe that $p(u)$ represents the cardinality of the set $S_1(u)$
of correlated columns with correlation higher than $\gamma/u$ (and
so parameterized by $u$), with associated characteristics $t^*(u)$,
$\gamma_{BT}(u)=\gamma/u$, $\gamma_{BG} =\gamma$. We are looking for
$u$  such that $$t^*(u)\gamma_{BT}(u)\le \gamma_{BG}
\Longleftrightarrow t^*(u)\gamma/u\le \gamma\Longleftrightarrow u\ge
u_1$$ since $t^*(u)=k/p(u)$ . Let us draw now the curve $p(u) \log
p(u)$ and find the point $u_2$ verifying
$$u_2=\inf \{u>0,\; p(u) \log p(u)\ge k,\; \}.$$
Deciding that the number of groups is
$$p^*=\lfloor u_1\vee u_2\rfloor ,$$
we are sure that the leading quantity in (\ref{tobalance}) is
$\gamma$ at least as soon as $\gamma\ge c\sqrt{\log p/n}$ which is
the standard case in high dimension.

\subsubsection{Determination of the delegates}
The set of 'delegates'
$${\mathcal D}=\left\{\ell \in\{1,\ldots,k\},\,\exists
\ell^\prime\in\{1,\ldots,k\}\setminus\{\ell\} \mbox{ such that } |
\Gamma_{\ell\ell^\prime}|>\gamma/p^* \right\}$$ is
also identified with the 'task' $t=1$. Each delegate is associated to one group.
It remains to distribute the variables whose indices are not in the ${\mathcal D}$
in the different $p^*$ groups.

\subsubsection{Completion of the groups}
 The variable of rank one in each group  ${\mathcal G}_j$ is a variable belonging
 to $\mathcal D$. The repartition is done in such a way that all the groups have
 the same (or almost the same) cardinality. In the same way as for the gathering
 Grouping procedures, we propose two versions for the Boosting Grouping:
\begin{itemize}
\item ({\bf BGc}): We rearrange the groups by sorting the correlation indicators
associated to the delegates:
 $R_{(1)}\geq \ldots\ge R_{(p^*)}$. This means that $\G_{1}$ contains the
 delegate $\ell_1$ such that $R_{\ell_1}=Y^t X_{.\ell_1}$ takes the largest
  correlation value (equal to $R_{(1)}$)
 and $\G_{p^*}$
has the delegate with the smallest $R_{\ell_{p^*}}$ correlation value (equal
to $R_{(p^*)}$).
The groups are then built such that the $R$'s are as homogeneous  as
possible in each group and as close as possible to their delegate. Grouping
starts by ranking the remaining $R$'s (i.e.
not  associated to a delegate): $R_{(1)} \geq \ldots \geq
R_{(k-p^*)}$. We denote $(\ell)$  the index associated to the
quantity $R_{(\ell)}$. The $p^*$ different groups are then
successively filled by using the ranking indices:
$$
\G_1=\{\,\ell_1,(1),\ldots,(\lfloor k/p^*\rfloor)-1\,\},\quad
\ldots\quad ,\G_{p^*}=\{\,\ell_{p^*},(k-p^*-\lfloor k/p^*\rfloor)+1),
\ldots,(k-p^*)\,\}.
$$

\item ({\bf BGa}):  It is the same procedure as ({\bf
BGc}) but using the $|R|$'s instead of the the $R$'s.
Notice that in this case, we have rearranged the groups by sorting the
absolute value of the correlation indicators
associated to the delegates.
\end{itemize}
Again, to
understand the improvement provided by the BGa and BGc in the next section,
we also consider
\begin{itemize}
\item ({\bf BGr}): the groups are filled up completed randomly.
The $k-p^*$ variables are spread out randomly into the $p^*$ groups.
\end{itemize}

\subsection{Quality of BRG}
Let us now consider the estimator $\hat\beta^*$ of $\beta$ obtained using the procedure GR-LOL combined with a pre-processing using BRG algorithm  to form the groups. Applying Theorem \ref{JAIME} under the conditions of the
theorem, it is easy to show that as soon as $\gamma\ge c[\log
p/n]^{1/2}$

\begin{align}\label{boost}
E\|\hat\beta^*-\beta\|_2^2&\leq
C\,\left(\gamma\right)^{2-q}.
\end{align}

\section{Simulation}\label{simu}
In this section, an extensive simulation study is conducted to
explore the practical qualities of procedure GR-LOL as well as the
Boosting Grouping (BRG) procedure. In the first part, we briefly describe  the experimental design and the empirical
tuning of the parameters of the procedure. The second part is
devoted to the study of the Boosting Grouping procedures comparing
to the gathered procedures given in Section \ref{Gather} and to the
procedures GGr and BGr where the groups are filled randomly.
Finally, GR-LOL procedure (with a pre BRG-processing) is compared
with two other procedures: LOL and Group Lasso. The comparison with
LOL (see \cite{MougeotPicardTribouley:2012}) allows to check the
contribution of the grouping and the comparison with the Group lasso
(see \cite{Yuan-Lin}) allows to evaluate GR-LOL with respect to this
 challenging procedure involving an important optimization step.

\subsection{Experimental design}
\subsubsection{Generation of the variables}
 The design matrix $X$ is a standard Gaussian $n \times
k$ matrix. Each column vector $X_{\cdot\ell}$ is centered and
normalized. The target observations $Y$ are given by $Y=X\beta+W$
where
\begin{itemize}
\item  $\beta$ is a vector of size $k$ whose coordinates are zero except $S$ which are
 $\beta_\ell = (-1)^{b_\ell} |z_\ell|$ for $\ell=1, \ldots, S$ where the
 $b$'s
 are i.i.d. Rademacher variables and the $z$'s are i.i.d. ${\cal N}(5,1)$
variables.
\item $W$ are i.i.d. variables ${\cal N}(0,\sigma^2)$. The
variance $\sigma^2$ of the noise is chosen such that the SNR
(signal over noise ratio) is close to $5$ which corresponds to a
middle noise level.
\end{itemize}

To introduce some dependency between the regressors, we choose
randomly a set denoted $\mathcal{R}$ of size $p_d=\lfloor\pi
k\rfloor$ of variables among the $k$ initial variables. Let us
denote by $M_{\rho}$ the $p_d \times p_d$ correlation matrix such
that $M_{\rho}(i,i)=1$ and $M_{\rho}(i,j)=\rho$ if $i\not=j$. Let
$V$ the eigenvector matrix and $D$ the diagonal eigenvalue matrix of
$M_{\rho}$ satisfying the singular value decomposition $M_{\rho}=V D
V^t$. Simulating a random gaussian matrix $Z$ of size $n \times
p_d$, we compute
 $X_{\mathcal R} = Z D^{1/2} V^t$; this resulting matrix has columns
 $X_\ell$ and $X_{\ell^\prime}$ verifying
 $cor(X_\ell,X_{\ell^\prime})=\rho$ as soon as
 $\ell \neq \ell^\prime$.  In order to study
broad experiments, different proportion values ($\pi=5\%, 10\%,
20\%$) as correlation values ($\rho=0.0,0.6,0.8$) have been studied.
This method has the advantage to tune accurately the number of
correlated variables as well as the amount of correlation between
the variables.

\subsubsection{Tuning parameters of the algorithms} As usual for thresholding
methods,  parameters $\lambda_n(1)$ and $\lambda_n(2)$ involved in
the GR-LOL procedure are critical values quite hard to tune because they
depend on constants which are not optimized and may not be available
 in practice. In this work, we tune them in an empirical way described  as follows:

\begin{description}
\item [Threshold $\lambda_n(1)$.] The first threshold  is used to select the leader
groups. Remember that at this stage, the number $p$ of groups is
known, (or determined by BG). Indeed, we do not determine directly
the level $\lambda_n(1)$ but find the number $p_0$ of leader groups
which is equivalent. Rearrange the groups along the values of
$\rho_j$ and denote $\G_{(1)},\ldots,\G_{(p)}$ the result of the
ranking. More precisely, the group $\G_{(j)}$ is associated to the
quantity $\rho_{(j)}$ where $\rho_{(j)}$ is the $j$th element of the
list $\rho_{(1)}^2 \geq \ldots \geq \rho_{(p)}^2$. We also denote
$t_{(j)}$ the cardinality of such a group $\G_{(j)}$ and  $p_0$ is
simply determined by
$$\sum_{j=1}^{p_0}t_{(j)}<n\quad\mbox{ and }\quad
\sum_{j=1}^{p_0+1}t_{(j)}\geq n.$$

When using grouping procedures, original variables are not handled
directly but through groups. If an important variable (i.e.
associated with a large coefficient of correlation with the target)
belongs to a cluster among unimportant variables (associated with
small coefficients), this variable may easily be unseen and killed
during the first thresholding step. This procedure slightly differs
from the LOL original procedure in being much less restrictive
during the first thresholding step and allowing to finally keep more
variables through the groups.

\item [Threshold $\lambda_n(2)$.] In order to compute the second thresholding step, we
do not determine, as previously, directly the level $\lambda_2(n)$
but find the number $p_1$ of finally retained groups which is
equivalent. The second threshold $\lambda_n(2)$ used for denoising
is computed by $5$-fold cross-validation. A proportion of $80\%$ of
the observations are used to estimate the $\beta$ coefficients.

The  $p_0$ groups, kept after the first thresholding, are ranked
using the $l^2$-norm of their estimated coefficients,
$\|\hat{\beta}\|_{\G_j,2}$. Each  $\G_{(j)}$ group, associated to
the quantity  $\|\hat{\beta}\|_{\G_j,2}$  is corresponding to the
$j$th element of the list $\|\hat{\beta}\|_{\G_(1),2} \geq \ldots
\geq \|\hat{\beta}\|_{\G_(p_0),2}$. The $20\%$ remaining
observations are used to sequentially compute the prediction error
using the one, ..., the $j$th first groups of the previous ranking
list. Using a model involving the $j$th first groups, the
 prediction error is defined by $\| Y - \hat{Y}_{{\cal U}_j}\|_2^2$
where ${\cal U}_j = \G_{(1)} \cup \ldots \cup \G_{(j)}$. The
prediction error is averaged using the $5$-fold cross-validation.
Finally, the first groups corresponding to the minimum prediction
error are kept.
\end{description}

In Section \ref{comparison}, we use  LOL  and the Group Lasso
algorithms which both tuning parameters as well. LOL algorithm is a
particularly case of GR-LOL when the number of groups equals the
number of variables i.e. $p=k$. For fair comparison, we use here for
LOL the same algorithm as for GR-LOL in the case where $p=k$. (And
so we have here a slight difference with the LOLA procedure provided
in \cite{MougeotPicardTribouley:2012}.)
 For group Lasso, the number of
final groups is computed by cross-validation as described in
(\cite{Yuan-Lin}, \cite{MaSong2007}, \cite{HuangHorowitz2010}). As
usual, the initial sample of observations is split into two samples:
the training set contains $75\%$ of the $n$ observations and is used
when the algorithms are running, the test set contains $25\%$ of the
$n$ observations and is used for the cross-validation methods.

\subsubsection{Criterion to evaluate the quality of the method}
For each studied procedure $P$ ($P$ is either $BG_{(a,c,r)}$ or
$GG_{(a,c,r)}$) with the prediction $\hat Y^P$, the
 relative prediction error $E_Y^P=\|Y - \hat {Y}^P\|_2^2/
\|Y\|_2^2$ is computed on the target $Y$. The results presented in
the tables give median values and standard deviations when  $K=100$
replications of the algorithms are performed. When GR-LOL is
compared with another procedure $P$ ($P$ is either LOL or Group
Lasso), the ratio $E_Y^P/E_Y^{GR-LOL}$ is computed. If the ratio is
close to $1$, the methods perform similarly;  when the ratio is
larger than $1$, GR-LOL outperforms P.

\subsubsection{BRG: Number $p$ of groups}
\begin{figure}[h]
\begin{center}
\includegraphics[width=8.5cm,height=4cm]{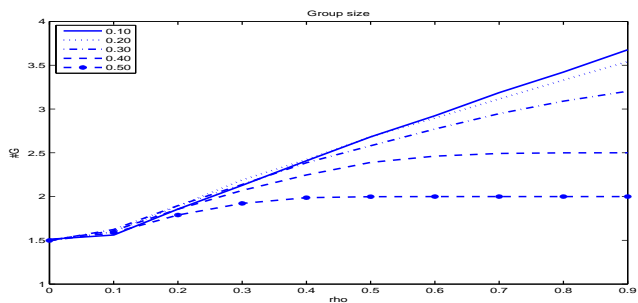}
\caption {$Y$-axis: Size of the groups. $X-$axis: correlation $\rho$
between the regressors for $\pi=10\%,20\%,30\%,40\%,50\%$. $n=200$,
$k=1000$, $SNR=5$, $K=100$. \label{groupsize}}
\end{center}
\end{figure}

Recall that the first step of  BRG  consists in  determining the
number $p^*$ of groups. and is detailed in Section
\ref{numberofgroups}. Figure \ref{groupsize} shows the average size
of the groups computed with the BRG procedure when the level of
dependence between the regressors given by $\pi$ and $\rho$ are
varying continuously.
 When no
dependency is introduced in the design matrix, we observe that the groups contain in
average $t^*=1.5$ variables using the experimental design previously described.
 Observe that the size of the groups is
increasing (and then the number $p^*$ of groups is decreasing) with
the level of dependency between the regressors (with $\pi$ or
$\rho$). For example, for $\rho=0.8$, the size of the group is
almost multiplied by 2 as $\pi$ decreases from $50\%$ to $10\%$.

For a fair comparison, the number $p^*$ of groups is the same for
all the methods, only the repartition of the variables between the
different groups varies. Defining the number of groups is not an
easy task. It should be underlined that in this case, the Random
Grouping and the Gathered Grouping both benefit of the optimal and
automatic choice of $p^*$ proposed by the boosting strategy. It should also be noticed
that the  Gathered and Boosting Grouping algorithms provide very different
configurations for the groups, as the average size $t^*$ of the groups
is small.
\subsubsection{Impact on the Coherence}

The empirical coherence $\gamma_{BG}$, $\gamma_{BT}$ and
$\gamma$ are computed and shown in Table \ref{cohgrouping} for
different value of correlation ($ \pi=0\%, 20\% , 40\%$ ; $\rho=0.0,
0.6, 0.8$) and for all the considered grouping strategies. For each
simulation, we have $\gamma=$sup$(\gamma_{BT},\gamma_{BG})$. As the
results presented in table \ref{cohgrouping} are averaged over
$K=100$ replications, we do not find necessarily at the end this
property, especially for Gathering grouping (GGa, GGc, GGr) which
can provide very different groups each time.

As expected, the boosting strategies induce a strong decrease of
$\gamma_{BT}$ as soon as there exists some dependency ($\pi > 20\%$,
$\rho > 0.6$). The different strategies for filling the groups (BGr,
BGc, BGa) does not have however any influence on $\gamma_{BG}$ as
expected also. The gathered groupings (GGc, GGa)  do not help to
reduce $\gamma_{BT}$ and $\tau^*$. As expected, the empirical
coherence (denoted $\gamma$ in the theoretical part) is increasing
with the dependence level $\rho$.
 Table \ref{cohgrouping} shows also the empirical value of
 $\tau^*=t^* \gamma_{BT} + \gamma_{BG}$ and
  $r^*=t^* \gamma_{BT}^2 + \gamma_{BG}^2$ computed for different strategies.

\begin{table}[h!]
\begin{center}
{\small
\begin{tabular}{ccccccc}
  \hline
$\pi=0\%$, $\rho=0.0$ & $t^*$ & $\gamma_{BT}$  & $\gamma_{BG}$ & $\gamma$ & $\tau^*$ & $r^*$ \\
  \hline
GGr  & 1.40 (0.06) & 0.321 (0.030)  & 0.317 (0.014) & 0.327 (0.015) & 0.766 & 0.245  \\
GGc  & 1.40 (0.06) & 0.318 (0.029)  & 0.321 (0.016) & 0.327 (0.015) & 0.766 & 0.245  \\
GGa  & 1.40 (0.06) & 0.318 (0.029)  & 0.319 (0.017) & 0.327 (0.015) & 0.764 & 0.243  \\
BGr  & 1.40 (0.06) & 0.234 (0.016)  & 0.327 (0.015) & 0.327 (0.015) & {\bf 0.655} & 0.184  \\
BGc  & 1.40 (0.06) & 0.234 (0.015)  & 0.327 (0.015) & 0.327 (0.015) & {\bf 0.655} & 0.184  \\
BGa  & 1.40 (0.06) & 0.234 (0.015)  & 0.327 (0.015) & 0.327 (0.015) & {\bf 0.655} & 0.184  \\
\hline \hline
$\pi=20\%$, $\rho=0.6$ & $t^*$ & $\gamma_{BT}$ & $\gamma_{BG}$ & $\gamma$ & $\tau$ & $\tau^*$ \\
\hline
GGr  & 2.80 (0.07) & 0.731 (0.029)  & 0.723 (0.020) & 0.733 (0.018) & 2.770 & 2.019  \\
GGc  & 2.80 (0.07) & 0.730 (0.032)  & 0.726 (0.019) & 0.733 (0.018) & 2.770 & 2.019  \\
GGa  & 2.80 (0.07) & 0.730 (0.032)  & 0.724 (0.019) & 0.733 (0.018) & 2.768 & 2.016  \\
BGr  & 2.80 (0.07) & {\bf 0.260} (0.018)  & 0.733 (0.018) & 0.733 (0.018) & {\bf 1.460} & 0.726  \\
BGc  & 2.80 (0.07) & {\bf 0.260} (0.016)  & 0.733 (0.018) & 0.733 (0.018) & {\bf 1.460} & 0.726  \\
BGa  & 2.80 (0.07) & {\bf 0.260} (0.016)  & 0.733 (0.018) & 0.733 (0.018) & {\bf 1.460} & 0.726  \\
\hline \hline
$\pi=40\%$, $\rho=0.6$ & $t^*$ & $\gamma_{BT}$  & $\gamma_{BG}$ & $\gamma$ & $\tau$ & $\tau^*$ \\
\hline
GGr  & 2.40 (0.02) & 0.742 (0.030)  & 0.739 (0.019) & 0.746 (0.019) & 2.521 & 1.869  \\
GGc  & 2.40 (0.02) & 0.742 (0.031)  & 0.740 (0.019) & 0.746 (0.019) & 2.522 & 1.870  \\
GGa  & 2.40 (0.02) & 0.741 (0.031)  & 0.740 (0.019) & 0.746 (0.019) & 2.519 & 1.866  \\
BGr  & 2.40 (0.02) & {\bf 0.303} (0.024)  & 0.746 (0.019) & 0.746 (0.019) & {\bf 1.474} & 0.778  \\
BGc  & 2.40 (0.02) & {\bf 0.303} (0.024)  & 0.746 (0.019) & 0.746 (0.019) & {\bf 1.474} & 0.778  \\
BGa  & 2.40 (0.02) & {\bf 0.303} (0.021)  & 0.746 (0.019) & 0.746 (0.019) & {\bf 1.473} & 0.777  \\
\hline \hline
$\pi=40\%$, $\rho=0.8$ & $t^*$ & $\gamma_{BT}$ &  $\gamma_{BG}$ & $\gamma$ & $\tau$ & $\tau^*$ \\
\hline
GGr  & 2.50 (0.03) & 0.868 (0.016)  & 0.866 (0.011) & 0.869 (0.011) & 3.035 & 2.631  \\
GGc  & 2.50 (0.03) & 0.867 (0.017)  & 0.867 (0.011) & 0.869 (0.011) & 3.035 & 2.632  \\
GGa  & 2.50 (0.03) & 0.867 (0.017)  & 0.867 (0.011) & 0.869 (0.011) & 3.035 & 2.632  \\
BGr  & 2.50 (0.03) & {\bf 0.315} (0.025)  & 0.869 (0.010) & 0.869 (0.010) & {\bf 1.656} & 1.003  \\
BGc  & 2.50 (0.03) & {\bf 0.314} (0.028)  & 0.869 (0.010) & 0.869 (0.010) & {\bf 1.654} & 1.002  \\
BGa  & 2.50 (0.03) & {\bf 0.317} (0.025)  & 0.869 (0.010) & 0.869 (0.010) & {\bf 1.662} & 1.007  \\
\end{tabular}
\caption {First line: Empirical coherence $\gamma_{BG}$, $\gamma_{BT}$, $\gamma$
 computed when the groups are built using the different strategies.
 $SNR=5$, $n=200$, $k=1000$, $K=100$,
$\pi=40\%$. \label{cohgrouping}} }
\end{center}
\end{table}

\subsubsection{Benefits of boosting grouping}
Table \ref{EpredBNRgrouping} compares the random (GGr), Gathered
(GGc, GGa) and boosting grouping (BGr, BGc, BGa) for different
sparsities $S$ and different levels of dependence ($\rho,\pi)$.
\begin{table}[h!]
\begin{center}
{\small
\begin{tabular}{cccccc}
  \hline
 $\pi=0\%,\rho=0$ & $S=10$ & $S=20$ & $S=30$ & $S=40$ & $S=50$ \\ \hline
GGr  &  3.06 ( 0.02) &  6.32 ( 0.04) & 11.36 ( 0.07) & 14.11 ( 0.09) & 14.53 ( 0.10)  \\
GGc  &  3.18 ( 0.02) &  5.18 ( 0.03) &  7.86 ( 0.04) &  9.81 ( 0.08) & 10.36 ( 0.06)  \\
GGa  &  2.97 ( 0.02) &  4.93 ( 0.02) &  7.62 ( 0.04) &  9.19 ( 0.07) & 10.31 ( 0.06)  \\
BGr  &  3.09 ( 0.02) &  7.09 ( 0.04) & 10.62 ( 0.06) & 11.64 ( 0.07) & 14.21 ( 0.09)  \\
BGc  &  3.15 ( 0.02) &  5.71 ( 0.05) &  9.84 ( 0.06) & 12.30 ( 0.07) & 12.31 ( 0.07)  \\
BGa  &  3.04 ( 0.02) &  5.33 ( 0.03) &  8.32 ( 0.05) & 10.12 ( 0.06) & 10.97 ( 0.07)  \\
 \hline \hline
 $\pi=20\%,\rho=0.6$ & $S=10$ & $S=20$ & $S=30$ & $S=40$ & $S=50$ \\ \hline
GGr  &  8.85 ( 0.16) & 28.82 ( 0.22) & 36.43 ( 0.23) & 40.61 ( 0.24) & 44.56 ( 0.23)  \\
GGc  &  7.52 ( 0.18) & 21.61 ( 0.24) & 41.72 ( 0.26) & 38.04 ( 0.26) & 37.50 ( 0.24)  \\
GGa  &  7.93 ( 0.17) & 24.43 ( 0.24) & 33.26 ( 0.26) & 40.67 ( 0.27) & 39.98 ( 0.24)  \\
BGr  &  9.35 ( 0.13) & 20.28 ( 0.17) & 25.35 ( 0.19) & 35.28 ( 0.20) & 32.80 ( 0.17)  \\
BGc  &  7.22 ( 0.10) & 14.41 ( 0.15) & 24.76 ( 0.16) & 25.50 ( 0.18) & 27.43 ( 0.19)  \\
BGa  &  {\bf 6.04 ( 0.05)} & {\bf 10.02 ( 0.06)} & {\bf 14.14 ( 0.10)} & {\bf 20.96 ( 0.13)} & {\bf 20.35 ( 0.14)}  \\
 \hline \hline
 $\pi=40\%,\rho=0.6$ & $S=10$ & $S=20$ & $S=30$ & $S=40$ & $S=50$ \\ \hline
GGr  & 19.78 ( 0.19) & 31.66 ( 0.20) & 39.23 ( 0.21) & 45.73 ( 0.22) & 44.66 ( 0.21)  \\
GGc  & 17.74 ( 0.18) & 38.77 ( 0.22) & 38.96 ( 0.23) & 51.72 ( 0.22) & 51.14 ( 0.23)  \\
GGa  & 18.28 ( 0.19) & 40.82 ( 0.22) & 43.42 ( 0.21) & 59.12 ( 0.22) & 54.80 ( 0.23)  \\
BGr  & 10.51 ( 0.07) & 17.34 ( 0.11) & 24.71 ( 0.13) & 30.16 ( 0.16) & 31.15 ( 0.19)  \\
BGc  &  9.48 ( 0.09) & 19.03 ( 0.14) & 24.26 ( 0.16) & 30.14 ( 0.17) & 31.27 ( 0.18)  \\
BGa  &  {\bf 7.51 ( 0.06)} &  {\bf 10.43 ( 0.07)} & {\bf 16.30 ( 0.09) }&  {\bf 20.63 ( 0.12)} &{\bf 23.41 (0.13) } \\
\hline \hline
 $\pi=40\%,\rho=0.8$ & $S=10$ & $S=20$ & $S=30$ & $S=40$ & $S=50$ \\ \hline
GGr  & 29.75 ( 0.20) & 43.95 ( 0.23) & 39.27 ( 0.23) & 48.75 ( 0.22) & 48.80 ( 0.27)  \\
GGc  & 37.59 ( 0.22) & 49.77 ( 0.25) & 49.81 ( 0.26) & 57.23 ( 0.24) & 53.57 ( 0.26)  \\
GGa  & 36.69 ( 0.21) & 51.53 ( 0.25) & 50.64 ( 0.26) & 59.99 ( 0.25) & 60.64 ( 0.26)  \\
BGr  & {\bf 7.85 ( 0.05)}&{\bf 13.95 ( 0.08)}&{\bf 18.14 ( 0.13)}&{\bf 19.82 ( 0.17)}&{\bf 26.48 ( 0.19)} \\
BGc  & {\bf 8.33 ( 0.07)}&{\bf 14.93 ( 0.11)}&{\bf 20.52 ( 0.15)}&{\bf 21.41 ( 0.17)}&{\bf 28.95 ( 0.18)} \\
BGa  & {\bf 5.96 ( 0.05)}&{\bf 9.19 ( 0.05)}& {\bf 12.72 ( 0.10)}&{\bf 16.26 ( 0.13)}&{\bf 19.44 ( 0.17)} \\
\hline
\end{tabular}
\caption {Relative prediction errors $E_Y$ ($\times 100$) for
Boosting Grouping (BGr, BGc, BGa), Gathered Grouping (GGc, GGa) and
Random Grouping (GGr) when the sparsity is varying, for various
levels of dependency given by $\pi,\rho$. $SNR=5$, $n=200$,
$k=1000$, $K=100$. \label{EpredBNRgrouping}} }
\end{center}
\end{table}
Let us first comment the no-dependency case ($\pi=0$). When the
sparsity is high ($S=10,20,30$), similar performances are obtained
for any grouping strategy. Underline that even building the groups
in a completely random manner is not a bad strategy. When the
sparsity is low ($S=40,50$), the Gathered Groupings (GGa and GGc)
bring the best results with a weak variability (low standard
deviation). As there is no specific correlation between the
regressors, the boosting procedure brings as expected in this case
no added value.

Actually, the boosting grouping procedure is especially adapted to
large correlation for taking advantage. For instance, when $\rho$
and $\pi$ are significative ($\rho=0.6$ and $\pi=0.4$), the boosting
procedure clearly shows substantial benefits. However, the
performances of the boosting depends on the strategy for filling the
groups. When the number of correlated variables is weak ($\pi=0.2$),
the boosting associated with  groups  filled randomly (BGr) is
rather competitive compared to Gathered groupings (GGc, GGa).
However, the boosting procedures with   groups filled homogeneously
always show  the best performances (BGc, BGa) with a preference for
the absolute value criteria. When there are strong correlations
between the regressors $\rho=0.6, 0.8$, the boosting procedures
(BRc, BGa) clearly outperforms  the random and the Gathered
grouping, and this is even true when the groups are filled  randomly
(BGR). BGa always brings the best results when the sparsity $S$
increases and/or the correlation $\rho$ between the regressors
increases.

\subsection{Study of the GR-LOL procedure}\label{comparison}

In this part, we present the performance results when  GR-LOL procedure  associated
with
 the Boosting Grouping strategy (BGa) is applied on the experimental design presented
 above. Comparisons between GR-LOL and LOL on the one hand, and GR-LOL and
Group-lasso on the second hand are explored.
\subsubsection{GR-LOL versus LOL}
The main difference between LOL and GR-LOL is that GR-LOL manipulates
groups of variables while LOL procedure handles  the variables directly.
 Table \ref{EpredGR-LOL-lol} shows a comparison of the performances obtained for LOL
  for the same experimental design
 as above.

\begin{table}[h!]
\begin{center}
{\small
\begin{tabular}{cc||ccccc}
  \hline
 $\pi$&$\rho$ &$S=10$ & $S=20$ & $S=30$ & $S=40$ & $S=50$ \\
 \hline \hline
$0\%$ & 0.00 & 1.083   & 1.522   & 1.616   & 1.777   & 1.666    \\
\hline
$20\%$ & 0.60 & 1.342   & 2.854   & 3.572   & 2.636   & 2.414    \\
$20\%$ & 0.80 & 1.877   & 5.436   & 3.898   & 3.117   & 2.649    \\
\hline
$40\%$ & 0.60 & 3.607   & 4.341   & 3.715   & 2.856   & 2.410    \\
$40\%$ & 0.80 & 6.287   & 6.429   & 4.773   & 3.440   & 3.417    \\

\hline
\end{tabular}
\caption {Relative prediction errors ratio $E_Y^{LOL}/E_Y^{GR-LOL}$
 for LOL and GR-LOL when the sparsity is varying for
different correlation values $\rho=0.0, 0.6, 0.8$ and rates $\pi=0,
0.2, 0.4$. $SNR=5$, $n=200$, $k=1000$. \label{EpredGR-LOL-lol}} }
\end{center}
\end{table}

We observe that LOL procedure performs particularly well when the
sparsity is large ($S$ small) and when the dependence between the
regressors is weak (\cite{MougeotPicardTribouley:2012}). In this
case, GR-LOL brings no improvement compared to LOL. Observe that, if
there is no dependency (case where $\rho=0.0$), the grouping
improves the performances of LOL when the sparsity decreases ($S$
increases). If the dependency increases (case where $\rho=0.6,0.8$),
GR-LOL always outperforms LOL for any considered sparsity.

\subsubsection{GR-LOL versus Group-lasso}
 The group Lasso
 is one of the most popular procedure for penalized regression
with grouping variables so we choose this method to challenge the
boosting Grouping procedure. To be fair, for both procedures, the
groups are built using the boosting strategy (BGa) and
cross-validation are both used to determine the final model.

Comparison of prediction results are given by Table \ref{EpredGR-LOLglasso}. Both
procedures show  similar behaviors in two cases:
 when there is no high correlation between the co variables
($\pi=0$) or when the sparsity ($S=50$) is small. In the other cases
(especially when the sparsity is large i.e. $S$ small), the results
given by GR-LOL are excellent: GR-LOL always outperforms the group
lasso.

\begin{table}[h!]
\begin{center}
{\small
\begin{tabular}{cc||ccccc}
 \hline
$\pi$&  $\rho$ & $S=10$ & $S=20$ & $S=30$ & $S=40$ & $S=50$ \\
 \hline
$0\%$ &0.0 & 1.228& 1.318 &1.143& 1.827 &2.001\\
\hline
$5\%$& 0.6 & 4.584 &2.366& 1.470& 1.944& 1.706\\
$5\%$& 0.8 &5.179& 2.490& 1.937 &1.122& 0.829\\
  \hline
$10\%$& 0.6 &2.764& 3.124& 1.892 &1.825 &0.967\\
$10\%$& 0.8 & 4.744& 1.643& 1.824& 1.511 &0.739\\
 \hline
$20\%$& 0.6  & 2.176& 3.032& 1.764& 1.385& 1.426\\
$20\%$& 0.8  & 3.250& 3.015& 1.986& 1.098& 1.048\\
 \hline
\end{tabular}
\caption {Relative prediction errors ratio $E_Y^{GLasso}/E_Y^{GR-LOL}$
for GR-LOL and GLasso when the sparsity is varying, for various levels
of dependency given by $\pi,\rho$. $SNR=5$, $n=200$, $k=1000$.
\label{EpredGR-LOLglasso}} }
\end{center}
\end{table}
To end this comparison, let us give a few words about computational
aspects. The Group Lasso algorithm is based on an optimization
procedure which can be time consuming while GR-LOL procedure solves
the penalized regression using two thresholding steps and a
classical regression. Regarding the complexity of the different
methods, GR-LOL has a strong advantage over the Group Lasso.

\subsection{Conclusion}

This experimental study shows that true  benefits can be obtained
using a grouping approach for penalized regression even in the case where there is
no prior knowledge on the groups. However, the results are highly relying on the
grouping strategy. The boosting strategy brings a nice
answer to the grouping problem when no prior
information is available on the structured sparsity. This strategy is very easy
to implement and
especially well adapted when a strong correlation exists between
the regressors in the case of high sparsity ($S$ small).

\section{Proofs}\label{sectionproof}

\subsection{RIP and associated properties:
$\tau^*$-conditions}\label{sectiontau}
In this part, we collect properties which are linked with the
coherence $\tau^*$. All these inequalities are extensively used in
the proof of Theorem \ref{JAIME} and the proofs of the propositions
stated in Section \ref{sectionr}; their proofs are detailed in the
appendix.

 Recall that for $\cI\subset\{1,\ldots,k\}$, $\Gamma_{\cI}=\tilde
X_\cI^t\tilde X_\cI$ is the associated Gram matrix of $\tilde
X_\cI$. $\tilde X_\cI$ is the matrix restricted to the columns of
$\tilde X$ whose indices are in $\cI$. Denote by $P_{V_\cI}$ the
projection on the space ${V_\cI}$ spanned by the predictors $\tilde
X_\ell$ whose indices $\ell$ belong to $\cI$. We also denote
$\bar\alpha(\cI)$ the vector of $\R^{\#(\cI)}$, such that
\begin{align}\label{baralpha}\tilde
X_\cI\bar\alpha(\cI)=P_{V_\cI}[\tilde X\alpha].\end{align} As well,
we define $\hat\alpha(\cI)$ the vector of $\R^{\#(\cI)}$, such that
\begin{align}\label{hatalpha}\tilde
X_\cI\hat\alpha(\cI)=P_{V_\cI}[Y].\end{align}

 The following lemma
describes the 'bloc-diagonal' aspect  of the Gram matrices
$\Gamma_\cI$ at least when the set of indices $\cI$ is small enough.
It is corresponding to the 'group-version' of the link between
coherence and RIP property (see for instance the corresponding result in
\cite{MougeotPicardTribouley:2012}).

\begin{lemma}\label{lemRIP} (RIP-property)
Let  $0<\nu<1$ be fixed. Let  $\cI$ be a subset of $\{1,\ldots,k\}$ such that $\tau(\cI)\leq \nu$.
 Then we get
\begin{equation}\label{cond}
\forall x\in \R^{\#(\cI)},\quad     \ \|x\|_{2}^2(1-\nu)\le x^t\,
\Gamma_{\cI}\,x \leq   \|x\|_{2}^2(1+\nu). \end{equation}
\end{lemma}
We deduce that the Gram
matrix $\Gamma_\cI$ is almost diagonal and in particular invertible as soon
as $\tau(\cI)\leq \nu$.
When this upper bound on $\tau(\cI)$ holds, we also extensively use the
RIP Property \eref{cond} in the
following forms:
\begin{equation}\label{condinv}
\forall x\in \R^{\#(\cI)},\; \|x\|_{2}^2(1+\nu)^{-1}\le
x^t\,\Gamma_\cI^{-1}\,x \leq  (1-\nu)^{-1} \|x\|_{2}^2 \,,
\end{equation}
and
\begin{equation}\label{rho-eucl}
\forall x\in \R^{\#(\cI)},\;
(1-\nu)\|x\|_{2}^2\le
\|\sum_{\ell\in\cI} x_\ell\;\tilde X_{\bullet\ell} \;\|_{2}^2\le
(1+\nu)\|x\|_{2}^2\,.
\end{equation}
 We also need the
following lemma
 \begin{lemma}\label{projobis}
For any $\cI$  subset of $\{1,\ldots,k\}$ such that $\tau(\cI)\leq \nu$, we have
 \begin{align}
 \forall x\in \R^{n},\;(1+\nu)^{-1}\sum_{\ell\,\in \cI}\;\left(\sum_{i=1}^n
  x_i\tilde X_{i,\ell}\right)^2\le
 \|P_{V_\cI}x\|_{2}^2\le (1-\nu)^{-1}\sum_{\ell\,\in
 \cI}\;\left(\sum_{i=1}^n
  x_i\tilde X_{i,\ell}\right)^2\,.
 \end{align}
 \end{lemma}

\subsection{Behavior of the projectors:
$r^*$-conditions}\label{sectionr}
In this subsection, we describe properties of the projection which
are more general as in the previous part where the results were
linked to the RIP property. These properties depend on the index
$r^*=t^*\;\gamma_{BT}^2+\gamma_{BG}^2$. It is noteworthy to observe
that in the no-group setting, we do not need to introduce this
 indicator $r^*$ since in this case $r^*=(\tau^*)^2$. Hence this
is one of the precise place where the grouping induces different
argument.

Let now state the following different technical results, which are
essential in the sequel.
\begin{lemma}\label{normB} Let $\cI,\cC$ be subsets of $\{1,\ldots,k\}$ and
 put
 $$ B(\cC)_\ell=\sum_{\ell^\prime\in\cC,\ell^\prime\not=\ell}
\;\Gamma_{\ell\ell^\prime}{\alpha}_{\ell^\prime} $$  for any $\ell$.
Then, we have
\begin{align*}
\|B(\cC)\|_{\cI,2}^2&\leq \;2\;\|\alpha\|_{\cC,1}^2\,r(\cI)
\end{align*}
where $r(\cI)$ is defined in (\ref{rI}).
\end{lemma}

\begin{proposition}\label{concR}For any integer $j$ from $\{1,\ldots,p\}$,
we get
\begin{align*}
\left|\, \| R\|_{\G_j,2}-\|\alpha\|_{\G_j,2}\right|^2\leq
4M^2v_nr(\G_j)+2(1+\nu)\|P_{V_{\G_j}}W\|_2^2.
\end{align*}
\end{proposition}

\begin{proposition}\label{projection} For any  subset $\cI$
of the leaders indices set $\G_\cB$, there exists $\kappa$ depending on $\nu$ such that
\begin{align*}
\|\widehat{\alpha} -{\alpha}\|_{\cI,2}^2\;\I\{\cI\subset
\G_\cB\}&\le \kappa\,\left(\|\alpha\|_1^2\,r(\cI) +
\|P_{V_{\cI}}W\|_{2}^2+\,\|P_{V_{\G_\cB}}W\|_{2}^2\,r(\cI)\right).
\end{align*}More precisely
$$
\kappa \geq \frac{1}{1-\nu}\vee \frac{6}{(1-\nu)^3}\vee
\frac{4(2\nu^2-\nu+2)}{(1-\nu)^4}.
$$
\end{proposition}
 \begin{proposition}\label{chi2}
 Let $\cI$ be a non random subset  such that
  $\#(\cI)\le n_{\cI}$, where $n_{\cI}$ is a deterministic quantity, then
   \begin{align}\label{expon}
 P\left(\frac{1}{\sigma^2}\|P_{V_{\cI}}[W]\|_{2}^2\geq
 z^2\right)&\leq \exp\left(-z^2/16\right)
 \end{align}
 for any $z$ such that $z^2\geq 4\,n_{\cI}$.
 If now $\cI$ is a random subset of the form
 $\{(j,t),\; j\in A,\;1\le t\le t_j\}$ where $A$ is a random set of
 $\{1,\ldots,p\}$ of cardinal less than $L$ (deterministic constant),
Inequality \eref{expon} is still true but
 for any $z$ such that $z^2\geq 16\,L\,(t^*\vee\log{p})$.
 In particular, this implies that for such a set,
 for any $k\ge 1$, there exists a constant $C_k$ such that
 \begin{align}\label{momk}
 E\left(\,\frac{1}{\sigma^2}\|P_{V_{\cI}}[W]\|_{2}^2\,\right)^k\le C_kL^k
 \,\left(t^*\vee\log{p}\right)^k.
 \end{align}
\end{proposition}

\subsection{Proof of the Theorem}
Thanks to Condition \eref{homo}, we have
\begin{align*}
a\nu_n\;\|\hat\beta^*-\beta\|_2^2\leq \|\hat\alpha^*-\alpha\|_2^2\leq b\nu_n\;\|\hat\beta^*-\beta\|_2^2
\end{align*}
which allows us to focus on the estimation error $\|\hat\alpha^*-\alpha\|_2$. We have
\begin{align*}
\|\hat\alpha^*-\alpha\|_2^2&=\|\hat\alpha^*-\alpha\|_{\G_\cB,2}^2+\|\alpha\|_{(\G_\cB)^c,2}^2
:= I \;(\mbox{In})\;+O \;(\mbox{Out}).
\end{align*}
We split $I$ into four terms :
\begin{align*}
I&= \sum_{j\in\cB}\I \{\|\hat\alpha\|_{\G_j,2}\ge \lam_n(2)\} \;
\I\{\|{\alpha}\|_{\G_j,2}\ge \lam_n(2)/2\}\;\|\hat\alpha-\alpha\|_{\G_j,2}^2 \\
& + \sum_{j\in\cB}\I \{\|\hat\alpha\|_{\G_j,2}\ge \lam_n(2)\}
\;\I\{\|{\alpha}\|_{\G_j,2}< \lam_n(2)/2\}
\;\|\hat\alpha-\alpha\|_{\G_j,2}^2 \\
&+\sum_{j\in\cB}\I\{\|\hat\alpha\|_{\G_j,2}<
\lam_n(2)\}\;\I\{\|\alpha\|_{\G_j,2}\ge
2\lam_n(2)\}\;\|\alpha\|_{\G_j,2}^2\\
&+\sum_{j\in\cB}\I\{\|\hat\alpha(\|_{\G_j,2}< \lam_n(2)\} \;
  \I\{\|\alpha\|_{\G_j,2}< 2\lam_n(2)\}  \; \|\alpha\|_{\G_j,2}^2\\
&:=IBB \;(\mbox{InBigBig})\;+\;IBS\;(\mbox{InBigSmall})\;+ISB
\;(\mbox{InSmallBig})\;+ISS
\;(\mbox{InSmallSmall})\;.
\end{align*}
We have on the other hand,
\begin{align*}
 O&\leq \sum_{j\in\cB^c}\I \{\|R\|_{\G_j,2}\le \lam_n(1)\}
\;\I\{\|\alpha\|_{\G_j,2}\ge 2\lam_n(1)\}\; \|\alpha\|_{\G_j,2}^2\\
&+\sum_{j\in\cB^c}\I \{\|R\|_{\G_j,2}\le \lam_n(1)\}
\;\I\{\|\alpha\|_{\G_j,2}< 2\lam_n(1)\}\; \|\alpha\|_{\G_j,2}^2\\
 &+\sum_{j\in\cB^c}\I \{\|R\|_{\G_j,2}\geq
\lam_n(1)\}
\;\I\{\|\alpha\|_{\G_j,2}\ge \lam_n(1)/2\}\; \|\alpha\|_{\G_j,2}^2\\
&+\sum_{j\in\cB^c}\I \{\|R\|_{\G_j,2}\geq \lam_n(1)\} \;
\I\{\|\alpha\|_{\G_j,2}< \lam_n(1)/2\}\; \|\alpha\|_{\G_j,2}^2\\
&:=OSB \;(\mbox{OutSmallBig})\;+OSS
\;(\mbox{OutSmallSmall})\;+OBB\;(\mbox{OutBigBig})\;+OBS
\;(\mbox{OutBigSmall})\;.
\end{align*}
\subsubsection{Study of $IBB$ and $ISB$}
Let us first study $ISB$. Observe that the two conditions
$\|\hat\alpha\|_{\G_j,2}\leq \lam_n(2)$ and $\|\alpha\|_{\G_j,2}\geq
2\lam_n(2)$ imply $\|\hat\alpha\|_{\G_j,2}\leq
\|\alpha\|_{\G_j,2}/2$. We deduce that
$$
\|\hat\alpha-\alpha\|_{\G_j,2}\geq
\|\alpha\|_{\G_j,2}-\|\hat\alpha\|_{\G_j,2}\geq
\|\alpha\|_{\G_j,2}/2
$$
and then
\begin{eqnarray}\label{ISB}
 ISB&\le& 4\,\sum_{j\in\cB}\I\{\|\alpha\|_{\G_j,2}\ge 2\lam_n(2)\}
  \;\|\hat\alpha-\alpha\|_{\G_j,2}^2= 4
  \;\|\hat\alpha-\alpha\|_{\cI\cap\G_\cB,2}^2
\end{eqnarray}
where \begin{align*}\cI:=\{(j,t)\in\{1,\ldots,k\},\;
\|\alpha\|_{\G_j,2}\ge 2\lam_n(2)\}.
\end{align*}
Thanks to Condition \ref{l1qbis}, we get
\begin{align*}
n_G(\cI):=\#(\,\{j,\; \exists t,\; (j,t)\in\cI\}\,)
&\leq\sum_{j=1}^p\I\{j\in\cC\},\;\left((2\lam_n(2))^{-1}\|\alpha\|_{\G_j,2}\right)^q
\leq (2\lam_n(2))^{-q}\,M^{q}\,\nu_n^{q/2}
\end{align*} and we bound $\#(\cI)$ by $t^*\times\#(\{j,\; \exists t,\;(j,t)\in\cI\})$.
It follows that
\begin{align*}
r(\cI)&\leq M^{q}\,v_n^{q/2}\,(2\lam_n(2))^{-q}\,[\gamma_{BG}^2+ t^*\gamma_{BT}^2]\\
 &\leq M^{q}\,v_n^{q/2}\,(2\lam_n(2))^{-q}\,r^*.
\end{align*}
Using successively Proposition \ref{projection} and Proposition \ref{chi2}, we get
\begin{align*}
 E(ISB)&\le 4\kappa E\left(M^2v_nr(\cI)+ \|P_{V_{\cI}}[W]\|_{2}^2+r(\cI)\|P_{V_{\G_\cB}}[W]\|_{2}^2\right)\\
  &\le 4\kappa \left([M^2v_nr(\cI)+C_1\sigma^2[t^*\vee \log p][n_G(\cI)+r(\cI)N^*]\right)
\\
 &\le [4+2C_1]\kappa\left(M^{ q}\,v_n^{q/2}\,(2\lam_n(2))^{-q}\right)
 \;(\lam^*)^2
\end{align*}
where $\lam^*$ is defined in (\ref{lam0}) and because $r^*N^*\leq\tau^*N^*\le
\nu$. The bound given in (\ref{ISB}) is valid for $IBB$ and then
the proof also holds for $IBB$.

\subsubsection{Study of $OSS$, $OBS$ and $ISS$}
Let $q$ be such that Condition (\ref{l1qbis}) is satisfied
\begin{align*}
OSS&\leq\sum_{j=1}^p\I\{\|\alpha\|_{\G_j,2}< 2\lam_n(1)\}\;
\|\alpha\|_{\G_j,2}^{2-q+q}\le (2\lam_n(1))^{2-q} \sum_{j=1}^p
\;\|\alpha\|_{\G_j,1}^{q}
\\&
\leq M^qv_n^{q/2}(2\lam_n(1))^{2-q}
\end{align*}
Note that this proof can also be  performed for $OBS$ and $ISS$
since $\lam_n(1)>\lam_n(2)$.
\subsubsection{Study of $OSB$}
Since
$$
\|\alpha\|_{\G_j,2}=\left(\|\alpha\|_{\G_j,2}-\|R\|_{\G_j,2}\right)+\|R\|_{\G_j,2}
$$
we get
\begin{align*}
OSB&\leq 2\sum_{j\in\cB^c}\I\{\|\alpha\|_{\G_j,2}-\|R\|_{\G_j,2}>
\lam_n(1)\}\; \;\left(\|\alpha\|_{\G_j,2}-\|R\|_{\G_j,2}\right)^{2}\\
&+2 \sum_{j\in\cB^c}\I\{\|\alpha\|_{\G_j,2}-\|R\|_{\G_j,2}>
\lam_n(1)\}\; \I\{\|R\|_{\G_j,2}\leq
\lam_n(1)\}\;\|R\|_{\G_j,2}^{2},
\end{align*}
and by Cauchy-Schwarz
\begin{align*}
E(OSB)&\leq
2\sum_{j\le p}\left[P\left(\,\left|\|\alpha\|_{\G_j,2}-\|R\|_{\G_j,2}\right|>
\lam_n(1)\right)\; \;E\left(\|\alpha\|_{\G_j,2}-\|R\|_{\G_j,2}\right)^{4}\right]^{1/2}\\
&+2 \lam_n(1)^2\sum_{j\le
p}P\left(\,\left|\|\alpha\|_{\G_j,2}-\|R\|_{\G_j,2}\right|>
\lam_n(1)\right).
\end{align*}
On the one hand, as an immediate consequence of Propositions \ref{concR} and \ref{chi2},
we get
\begin{align*}
E\left(\;\, \| R\|_{\G_j,2}-\|\alpha\|_{\G_j,2}\right)^4&\le
32M^4v_n^2\tau(\G_j)^2+16\sigma^4(1+\nu)^2\#(\G_j)^2.
\end{align*}
Since $\#(\G_j)\leq t^*$ and $\tau(\G_j)\leq\tau^*$, we bound this
term by $32(\lam^*)^4$. On the other hand, using Proposition
\ref{chi2},
\begin{align*}P\left(\;\left|\,\|R\|_{\G_j,2}|-\|\alpha\|_{\G_j,2}\,\right|\ge
\lam\right)
&\le P\left(\|P_{V_{\G_j}} W\|_{2}\ge \lam/2(1+\nu)^{1/2}\right)\\
&\leq  \exp \left(- \lam^2/32(1+\nu)\right)
\end{align*}
 as soon as
 $$\lam^2\geq \left(8 M^2v_n\,r(\G_j)\right)\vee
 \left(16(1+\nu)\,[t^*\vee \log p]\right).$$ This condition is verified by
  $\lam_n(1)$ as soon as \begin{align}\label{lam1}\lam_n(1)\geq 4\lam^*\end{align}
and it follows
\begin{align*}
E(OSB)&\leq \frac
94\;p\lam_n(1)^2\;\exp\left(-\frac{\lam_n(1)^2}{64(1+\nu)}\right).
\end{align*}
 \subsubsection{Study of $OBB$}
Let us decompose again this term into 2 different ones,
$$OBB=\sum_{j\in\cC_1}\|\alpha\|^2_{\G_j,2}+\sum_{j\in\cC_2}\|\alpha\|^2_{\G_j,2}
:=OBB_1+OBB_2$$ where
$$\cC=\{j\in \cB^c,\;\|\alpha\|_{\G_j,2}\ge 2\lam_n(1),\, \|R\|_{\G_j,2}\ge
\lam_n(1)\}$$ and
$$\cC_1=\cC\cap \{j, \|R\|_{\G_j,2}\le\|\alpha\|_{\G_j,2}/2\}\quad
\mbox{ and }\quad \cC_2=\cC\cap \{j,
\|R\|_{\G_j,2}\geq\|\alpha\|_{\G_j,2}/2\}
$$
On the one hand, we obviously have
\begin{align*}
\cC_1&\subset \{j\in \{1,\ldots,p\},\; \lam_n(1)\le
\|\alpha\|_{\G_j,2}-\| R\|_{\G_j,2}\}
\end{align*}
leading to
\begin{align*}
OBB_1&\leq \sum_{j=1}^p\I\{\left|\,\|\alpha\|_{\G_j,2}-\|
R\|_{\G_j,2}\right|\geq \lam_n(1)\}\;\|\alpha\|^2_{\G_j,2}.
\end{align*}
We conclude as for the term $OSB$. For $OBB_1$ the argument is
slightly more subtle: on the other hand,
$$j\not\in \cB\quad \mbox{ and }\quad \|R\|_{\G_j,2}\geq \lam_n(1)
\Longrightarrow \rho_{(j)}^2\leq \rho_{(N^*)}^2
$$
(see Step 2 of the procedure) inducing that there exist at least
$N^*$ leader indices $j^\prime\not= j$ in $\{1,\ldots,p\}$ such that
$\|R\|_{\G_{j^\prime},2}\geq \|R\|_{\G_j,2}$. Assume now that the following inequality
is true (this will be proved later):
\begin{align}\label{hypo_a_verifier}
\#(\cC)<N^*.
\end{align}
This implies that there exists at least one index (depending on $j$)
called $j^*(j)$ such that
$$
 \|\alpha\|_{\G_{j^*(j)},2}< \lam_n(1)/2 \;\;(\mbox{because } j^*(j)\not\in\cC)\;\; \quad \mbox{
and } \quad \| R\|_{\G_{j^*(j)},2}\ge \|R\|_{\G_j,2} .
$$
We deduce that, for this index $j^*(j)$, we have
\begin{align*}
\| R\|_{\G_{j^*(j)},2}-\| \alpha\|_{\G_{j^*(j)},2}&>\|
R\|_{\G_{j},2}- \lam_n(1)/2 \\
&>\| R\|_{\G_{j},2}/2\hspace{2.4cm}(\mbox{because }j\in\cC)\\
&>\| \alpha\|_{\G_{j},2}/2- \lam_n(1)/2\quad \quad\,(\mbox{because
}j\in\cC_2)\\
&>\| \alpha\|_{\G_{j},2}/4 \hspace{2.4cm}(\mbox{because }j\in\cC).
\end{align*}
It follows that
\begin{align*}
OBB_2&\leq 4\sum_{j=1}^p\I\{\left|\,\|\alpha\|_{\G_{j^*(j)},2}-\|
R\|_{\G_{j^*(j)},2}\right|\geq \lam_n(1)/2\}\;
\left(\,\|\alpha\|_{\G_{j^*(j)},2}-\| R\|_{\G_{j^*(j)},2}\right)^2
\end{align*}
and we conclude as for the term $OSB$. It remains now to prove
\eref{hypo_a_verifier}: thanks to Condition \ref{l1qbis}, we get
\begin{align*}
\#(\cC)&\leq \#(\,\{j\in \{1,\ldots,p\},\;\|\alpha\|_{\G_j,2}\ge
2\lam_n(1)\}\,)\\
&\leq\sum_{j=1}^p\I\{j\in\cC\},\;\left(2\lam_n(1)^{-1}\|\alpha\|_{\G_j,2}\right)^q
\leq \left(2\lam_n(1)^{-1}\right)^q\,M^{ q}\,\nu_n^{q/2}
\end{align*}
and  (\ref{hypo_a_verifier}) is satisfied as soon as $\lam_n(1)\geq
2\,M^{ }\,v_n^{1/2}\,(N^*)^{-1/q}$
which is verified for any $q\le 1$ as soon as
\begin{align}\label{lam2}
\lam_n(1)\geq 2\,M^{ }\,v_n^{1/2}\,(N^*)^{-1}.
\end{align}
\subsubsection{Study of $IBS$}
The triangular inequality for the norm $\|.\|_{\G_j,2}$ leads to
\begin{align*}
IBS&\leq \sum_{j\in\cB} \I \{\|\hat\alpha-\alpha\|_{\G_j,2}\ge
\lam_n(2)/2\}\; \|\hat\alpha-\alpha\|_{\G_j,2}^2.
\end{align*}
Using Cauchy Schwarz inequality
we get
\begin{align*}
E(IBS)&\leq \sum_{j=1}^p
\left(E\|\hat\alpha-\alpha\|_{\G_j,2}^4\,\I\{j\in\cB\}\right)^{1/2}\;
P\left(\,\|\hat\alpha-\alpha\|_{\G_j,2}\,\I\{j\in\cB\}\ge
\lam_n(2)/2\right)^{1/2}.
\end{align*}
On the one hand, by Propositions \ref{projection} and \ref{chi2}, we get
\begin{align*}
 E\left(\|\hat\alpha-\alpha\|_{\G_j,2}^4\,\I\{j\in\cB\}\right)&\le 3
 \kappa^2\left(M^4v_n^2r(\G_j)^2
 +2C_2\sigma^4\left[1+r(\G_j)^2[N^*]^2\right][t^*\vee \log p]^2\right).
\end{align*}
Since $r(\G_j)\leq r^*$, $\#(\G_j)\leq t^*$ and $\#(\G_\cB)\leq N^*t^*$, we bound this
term by $64\kappa^2\,(\lam^*)^4$ (using  $N^*\tau^*<\nu$). On the other hand, using again  Propositions \ref{projection} and \ref{chi2}, we have
\begin{align*}
 P\left(\;\|\alpha-\hat\alpha\|_{\G_j,2}\;\I\{j\in
\cB\}\,\ge\lam\;\right)&\le
\exp\left(-\frac{\lam^2}{24\kappa}\right)+
\exp\left(-\frac{\lam^2}{24\kappa r(\G_j)}\right)
\end{align*}
as soon as
$$
\lam^2\geq 3\kappa\left(M^2v_nr(\G_j)\vee
16(1+r(\G_j)N^*)[t^*\vee\log p]\right).
$$
It follows that, if
\begin{align}\label{lam3}\lam_n(2)\geq 5\sqrt{\kappa}\;\lam^*\end{align}
we get
\begin{align*}
E(IBS)&\leq 3\kappa\;p
(\lam^*)^2\left[\exp\left(-\frac{\lam^2(2)}{192\kappa}\right)+
\exp\left(-\frac{\lam^2(2)}{192\kappa r^*}\right)\right].
\end{align*}

\subsubsection{End of the proof}
If we summarize the results obtained above, choosing
$$\lam_n(1)=c_1\lam^*\vee [2\,M^{ }\,v_n^{1/2}\,(N^*)^{-1}]\quad\mbox{
and }\quad \lam_n(2)=c_2\lam^*$$ with $c_1>c_2$, $c_2\geq
5\sqrt{\kappa}$ and $ c_1\geq (4+\nu^{-1/q})$,
 we obtain
\begin{align*}
E\|\hat\alpha^*-\alpha\|_2^2&\leq 2\times \left[4\kappa\left(M^{
q}\,v_n^{q/2}\,(2\lam_n(2))^{-q}\right)
 \;(\lam^*)^2\right]+3\times\left[M^qv_n^{q/2}(2\lam_n(1))^{2-q}\right]\\
 &+\left[\frac 94\;p\lam_n(1)^2\;\exp\left(-\frac{\lam_n(1)^2}{64(1+\nu)}\right)\right]\\
&+\left[
3\kappa\;p
(\lam^*)^2\left[\exp\left(-\frac{\lam^2(2)}{64\kappa}\right)+
\exp\left(-\frac{\lam^2(2)}{192\kappa\tau^*}\right)\right]\right]\\
&\leq cv_n\left((\lambda^*)^{-q}v_n^{q/2-1}+(N^*)^{q-2}\right)
\end{align*}
under the condition
\begin{align*}
c_a\;p\;\exp\left(-c_b(\lam^*)(1\wedge (r^*)^{-1})\right)&\leq
v_n^{q/2}\,(\lam^*)^{-q}
\end{align*}
where
$$c_a=M^{-q}\left(\frac 94 c_1^2\vee 3\kappa\right)\mbox{ and }
c_b=\frac{c_1^2}{64(1+\nu)}\wedge \frac{c_2^2}{192\kappa}. $$
Replacing $\lambda^*$, we obtain the announced result.

\section{Appendix}

 Recall that $\tilde
X_\cI$ denotes the matrix restricted to the columns of $\tilde X$
whose indices are in $\cI$ subset of $\{1,\ldots,k\}$ and that $\Gamma_{\cI}=\tilde
X_\cI^t\tilde
X_\cI$. Denote
$P_{V_\cI}$ the projection on the space spanned by the predictors
$\tilde X_\ell$ whose index $\ell$ belongs to $\cI$
$$P_{V_\cI}=\tilde X_\cI(\tilde X_\cI^t\tilde X_\cI)^{-1}\tilde
X_\cI^t=\tilde X_\cI(\Gamma_\cI)^{-1}\tilde X_\cI^t.
$$
Recall that any index $\ell$ of $\{1,\ldots,k\}$ can be registered
as a pair $(j,t)$ where $j$ is the index of the group $\G_j$ where
$\ell$ is belonging and $t$ is the rank of $\ell$ inside $\G_j$.

\subsection{Proof of Lemma \ref{normB}}
 We use the definitions (\ref{btask}) and (\ref{bgroup}) of $\gamma_{BT}$ and $\gamma_{BG}$
\begin{align*}
\|B(\cC)\|_{\cI}^2&=\sum_{\ell\in\;
\cI}B(\cC)_\ell^2=\sum_{\ell\in\;
\cI}\left(\sum_{\ell^\prime\in\cC,\ell^\prime\not=\ell}
\;\Gamma_{\ell\ell^\prime}{\alpha}_{\ell^\prime} \right)^2\\&\leq \;\sum_{(j,t)\in\cI} \left(
\gamma_{BT}\;\sum_{(j^\prime,t^\prime)\in\cC,t^\prime\not
=t}|{\alpha}_{(j^\prime,t^\prime)}|
+\gamma_{BG}\;\sum_{j^\prime=1,\ldots,p,(j^\prime,t)\in\cC,j^\prime\not
=j}|{\alpha}_{(j^\prime,t)}|
\right)^2\\
&\leq \;2\,\gamma_{BT}^2\;\sum_{(j,t)\in\cI} \;\left(\sum_{\ell
\in\cC} |\alpha_{\ell}|\right)^2
+2\gamma_{BG}^2\;\sum_{j=1,\ldots,p,(j,t)\in\cI}
\left[\sum_{t=1,\ldots,t_j,(j,t)\in\cI}
 \;\left( \sum_{j^\prime=1,\ldots,p,(j^\prime,t)\in\cC,j^\prime\not
=j}|{\alpha}_{(j^\prime,t)}|\right)^2\right]
\\
&\leq \;2\,\gamma_{BT}^2\;\#(\cI) \;\left(\sum_{\ell
\in\cC}
|\alpha_{\ell}|\right)^2
+2\gamma_{BG}^2\;\#(\{j,(j,t)\in\cI\})\left(\sum_{\ell
\in\cC}
|\alpha_{\ell}|\right)^2
\end{align*}
which ends the proof.
\subsection{Proof of Lemma \ref{lemRIP}}
 Let us decompose the sum
\begin{align*}
x^t\,\Gamma_{\cI}\,x&=\sum_{\ell,\ell^\prime=1}^mx_{\ell}x_{\ell^\prime}
\left(\Gamma_{\cI}\right)_{\ell\ell^\prime}= \sum_{\ell=1}^mx_{\ell}^2
\left(\Gamma_{\cI}\right)_{\ell\ell}+\sum_{\ell\not=\ell^\prime}x_{\ell}x_{\ell^\prime}
\left(\Gamma_{\cI}\right)_{\ell\ell^\prime}.
\end{align*}
Using Condition \eref{norm}, it follows that
\begin{align*}
|x^t\,\Gamma_{\cI}\,x-\|x\|_{2}^2|&=
\sum_{\ell,\ell^\prime=1,\ldots,m,\ell\not=\ell^\prime}x_{\ell}x_{\ell^\prime}
\left(\Gamma_{\cI}\right)_{\ell\ell^\prime}.
\end{align*}
In order to solve the difficulty due to the fact that the size $t_j$
of the groups $\G_j$ could be different, we consider that $t$ is
varying until $t^*=\max(t_1,\ldots,t_p)$ with the convention that
$x_{(j,t)}=0$ if the index $(j,t)\not \in \cI$. Using Definition
(\ref{bgroup}) and Definition (\ref{btask}), we get
\begin{align*}
|x^t\Gamma(\cI)x-\|x\|_{{ l}_2(m)}^2| &\leq
\gamma_{BT}\;\sum_{(j,t)\in\cI,(j^\prime,t^\prime)\in\cI,t\not
=t^\prime}|x_{(j^\prime,t)}x_{(j^\prime,t^\prime)}|+\gamma_{BG}\;\sum_{(j,t)\in\cI,(j^\prime,t^\prime)\in\cI,t
=t^\prime }|x_{(j,t)}x_{(j^\prime,t)}|\\
&\leq \gamma_{BT}\; \left(\sum_{(j,t)\in\cI}|x_{(j,t)}|\right)^2+\gamma_{BG}\;\sum_{t=0}^{t^*}\left(\sum_{j\in \{j,(j,t)\in\cI\}}|x_{(j,t)}|\right)^2\\
& \leq \gamma_{BT}\;\#(\cI)\sum_{(j,t)\in\cI}|x_{(j,t)}|^2+\gamma_{BG}\;
\sum_{t=0}^{t^*}\,\#(\{j,(j,t)\in\cI\})\; \sum_{j\in\{j,(j,t)\in\cI\}} |x_{(j,t)}|^2\\
&\leq \tau(\cI) \;\|x\|_{2}^2
\end{align*}
which ends the proof since $\tau(\cI)\leq \nu$.

\subsection{Proof of Lemma\ref{projobis}}
Since
\begin{align*}
\|P_{V_{\cI}}x\|_{2}^2 &=( \tilde X_{\cI}^t x)^t\;(
\Gamma_\cI)^{-1}\;( \tilde X_{\cI}^t x),
\end{align*}
we have
\begin{eqnarray*}
 (1+\nu)^{-1}\,\|\tilde X_{\cI}^t x\|_{2}^2\leq
\|P_{V_{\cI}}x\|_{2}^2 &\leq (1-\nu)^{-1}\,\|\tilde X_{\cI}^t
x\|_{2}^2
\end{eqnarray*}
applying the RIP Property \eref{condinv}. Observing that
\begin{eqnarray*}
\|\tilde X_{\cI}^t x\|_{2}^2=( \tilde X_{\cI}^t x)^t\,(
\tilde X_{\cI}^t x)= \sum_{\ell\,\in \cI}\left(\sum_{i=1}^nx_i\tilde
X_{i,\ell }\right)^2,
\end{eqnarray*}
we obtain the announced result.

\subsection{Proof of Proposition \ref{concR}}
Since the model under consideration is $Y=\tilde X \alpha+W$, we have for any $\ell$ in $\{1,\ldots,k\}$
$$ R_\ell=\sum_{i=1}^n Y_i\tilde X_{i,\ell}=\sum_{i=1}^n (\tilde X_i\alpha)\tilde X_{i,\ell}+\sum_{i=1}^n W_i\tilde X_{i,\ell}
$$ leading to
\begin{align*}
 R_\ell-\alpha_\ell&=\sum_{\ell^\prime=1,\ldots,k,\ell^\prime\not=\ell}
 \Gamma_{\ell\ell^\prime}{\alpha}_{\ell^\prime}+\tilde X_{\ell}^t\,W := B_\ell+V_\ell
\end{align*}
thanks to Condition (\ref{norm}). It follows
\begin{align*}
\left|\|R\|_{\G_j,2}-\|\alpha\|_{\G_j,2}\,\right|&\leq
\|R-\alpha\|_{\G_j,2}\\
&\leq \|B\|_{\G_j,2}+\|V\|_{\G_j,2}.
\end{align*}
Applying Lemma \ref{normB} with $\cC=\{1,\ldots,k\}$ and $\cI=\G_j$,
we obtain $\|B\|_{\G_j,2}^2\leq \;2\,\|\alpha\|_{1}^2\,r(\G_j)$.
Since
\begin{equation}\label{l1ql1}
\|\alpha\|_{1}=\sum_{j=1}^p\sum_{t=0}^{t_j}
|{\alpha}_{(j,t)}|\le\left[\sum_{j=1}^p\|\alpha\|_{\G_j,1}^{q}\right]^{1/q}\le
Mv_n^{1/2},
\end{equation}
we get
\begin{align*}
\|B\|_{\G_j,2}^2&\leq  2M^{2}v_n\,r(\G_j).
\end{align*}
Second, using Property \eref{condinv} which holds because we assumed
\eref{mini}, we get
\begin{align*}
\|V\|_{\G_j,2}^2&=
\sum_{\ell\in\G_j}\left(\tilde X_{\ell}^t\,W\right)^2=
\sum_{\ell\in\G_j}\left(\tilde
X_{\G_j}^tW\right)_{\ell}^2
\\
&\leq (1+\nu)\,\left(\tilde
X_{\G_j}^tW\right)^t\Gamma_{\G_j}^{-1}\left(\tilde
X_{\G_j}^tW\right) =(1+\nu)\|P_{V_{\G_j}} W\|_{2}^2
\end{align*}
which ends the proof.

\subsection{Proof of Proposition \ref{projection}}
Recall the definitions \eref{baralpha} and \eref{hatalpha} and let us put
$$
\bar{\alpha}(\cI)={\alpha}_{\cI}+(\tilde X_{\cI}^t\tilde
X_{\cI})^{-1}\tilde X_{\cI}\;\tilde X_{\cI^c}{\alpha}_{\cI^c}
$$such that
$$
\bar{\alpha}(\cI)-{\alpha}_{\cI}=(\Gamma_{\cI})^{-1}\tilde
X_{\cI}\;\tilde X_{\cI^c}{\alpha}_{\cI^c}.
$$
Since $\cI\subset \G_\cB$, we have $\widehat{\alpha}_\ell=\widehat{\alpha}
(\cB)_{\ell}$ for any $\ell\in\cI$ and
\begin{align*}
\|{\alpha}
-\widehat{\alpha}\|_{\cI,2}^2&
=\|{\alpha}_{\cI}-\widehat{\alpha}(\cB)_{\cI}\|^2_{\cI,2}\\
& \leq \|\alpha_{\cI}-\bar{\alpha}(\cI)\|^2_{\cI,2}+
\|\bar{\alpha}(\cI)-\widehat{\alpha}(\cI)\|^2_{\cI,2} +
\|\widehat{\alpha}(\cI)-\widehat{{\alpha}(\cB)}_{\cI}\|^2_{\cI,2}\\
&:=t_1(\cI)+t_2(\cI)+t_3.
\end{align*}
Using twice the RIP Property and applying Lemma \ref{normB}
for $\cI:=\cI$ and $\cC:=\cI^c$, we bound the first term
\begin{align*}
t_1(\cI)&\leq \frac{1}{1-\nu}\;(\bar{{\alpha}}(\cI)-{\alpha}_{\cI})^t
\,\Gamma_{\cI}\,(\bar{\alpha}(\cI)-{\alpha}_{\cI})\\
&=\frac{1}{1-\nu}\;\;({\alpha}_{\cI^c}^t\tilde
X_{\cI^c}^t\tilde X_{\cI}^t) \,(\Gamma_\cI)^{-1}\,(\tilde
X_{\cI}\tilde X_{\cI^c}
{\alpha}_{\cI^c})\\
&\leq \frac{1}{(1-\nu)^2}\|\tilde
X_{\cI}^t\,\tilde X_{\cI^c}{\alpha}_{\cI^c}
\|_{\cI,2}^2\\
&= \frac{1}{(1-\nu)^2}\sum_{\ell\in
\cI}\left(\sum_{\ell^\prime\in \cI^c} \Gamma_{\ell\ell^\prime}\,{\alpha}_{\ell^\prime}\right)^2\\
&\leq
\frac{2}{(1-\nu)^2}\;\|\alpha\|_{\cI^c,1}^2
\,r(\cI).
\end{align*}
Recall that in the specific case where $\cI=\G_\cB$, we
get $\,\tau(\G_\cB)\leq \nu$ by construction of
the leader groups (see (\ref{leader})), so that
\begin{align}\label{t1B}
t_1(\G_\cB)&\leq
\frac{2\nu}{(1-\nu)^2}\;\|\alpha\|_{\G_\cB^c,1}^2.
\end{align}
For the study of $t_2(\cI)$, use Inequality (\ref{rho-eucl})
\begin{align*}
t_2(\cI)
 &\leq
\frac{1}{1-\nu}\;\|\tilde X_{\cI}\bar{\alpha}(\cI)-\tilde X_{\cI}\,\widehat{\alpha}(\cI)\|^2_2
\end{align*} and observe that
\begin{align*}
\tilde X_{\cI}\,\hat\alpha(\cI)
 =P_{V_{\cI}}[\tilde X\alpha+W]=\tilde X_{\cI} \bar{\alpha}(\cI)+P_{V_{\cI}}W
\end{align*}
to obtain the bound
\begin{align}\label{t2}
t_2(\cI)
&\leq \frac{1}{1-\nu}\;\|P_{V_{\cI}}W\|_{2}^2.
\end{align}
 finally, use again Inequality (\ref{rho-eucl})
\begin{align*}
t_3&\leq \frac{1}{1-\nu}\;\|\tilde X_{\cI}\widehat{\alpha}(\cI)-
\tilde X_{\cI}\widehat{\alpha}(\cB)_{\cI}\|^2_{2}\end{align*}
and observe that
\begin{align*}
\tilde X_{\cI}\widehat{\alpha}(\cI)-
\tilde X_{\cI}\widehat{\alpha}(\cB)_{\cI}&=P_{V_\cI}[\tilde X_{\cI}\widehat{\alpha}(\cI)-
\tilde X_{\cI}\widehat{\alpha}(\cB)_{\cI}]\\
&=P_{V_\cI}[\tilde X_{\cI}\widehat{\alpha}(\cI)-
\tilde X_{\G_\cB}\widehat{\alpha}(\cB)+\tilde X_{\G_\cB/\cI}\widehat{\alpha(\cB)}_{\G_\cB/\cI}]\\
&=P_{V_\cI}[P_{V_\cI}[\tilde X{\alpha}+W]- P_{V_{\G_\cB}}[\tilde
X{\alpha}+W]
+\tilde X_{\G_\cB/\cI}\widehat{\alpha}(\cB)_{\G_\cB/\cI}]\\
&=P_{V_\cI}[\tilde
X_{\G_\cB/\cI}\widehat{\alpha}(\cB)_{\G_\cB/\cI}].
\end{align*}
since $\cI\subset \G_\cB$. Applying Lemma \ref{projobis} and Lemma \ref{normB} for
$\cI:=\cI$ and $\cC:=\G_\B/\cI$, we get
\begin{align*}
t_3&\leq \frac{1}{1-\nu}\;\| P_{V_\cI}[\tilde
X_{\G_\cB/\cI}\widehat{\alpha}(\cB)_{\G_\cB/\cI}] \|^2_{2}
\\
&\le
\frac{1}{(1-\nu)^2}\;\sum_{\ell\in \cI}
\left(\sum_{i=1}^n\left(\sum_{\ell^\prime\in
\G_\cB/\cI}(\widehat{\alpha}(\cB))_{\ell^\prime}\tilde
X_{i\ell^\prime}\right)
\tilde X_{i\ell}\right)^2\\
&\leq \frac{2}{(1-\nu)^2}\; r(\cI)
\;\|\widehat{\alpha}(\cB)\|_{\G_\cB,1}^2.
\end{align*}
Writing
$$
\|\widehat{\alpha}(\cB)\|_{\G_\cB,1}^2\leq
3\left(\;\|\widehat{\alpha}(\cB)-\bar{\alpha}(\cB)\|_{\G_\cB,1}^2+
\|\bar{\alpha}(\cB)-\alpha\|_{\G_\cB,1}^2+\|\alpha\|_{\G_\cB,1}^2\;\right)
$$
we deduce that
\begin{align*}
t_3&\leq \frac{6}{(1-\nu)^2}\;r(\cI)\,
\left(t_2(\G_\cB)+t_1(\G_\cB)+\|\alpha\|_{\G_\cB,1}^2\right)
\end{align*}
and combining with (\ref{t1B}) and (\ref{t2}), we obtain
\begin{align*}
t_3 &\leq \frac{6}{(1-\nu)^2}\;r(\cI)
\,\left(\frac{2\nu}{(1-\nu)^2}\;
\|\alpha\|_{\G_\cB^c,1}^2+ \frac{1}{1-\nu}\|P_{V_{\G_\cB}}W\|_{2}^2
+\|\alpha\|_{\G_\cB,1}^2\right).
\end{align*}
This ends the proof of the proposition.

\subsection{Proof of Proposition \protect{\ref{chi2}}}
 first, the proof concerning the case where $\cI$ is not random is standard and can be found for instance
 in \cite{MougeotPicardTribouley:2012}.
Assume now that $\cI$ is random. We take into account all the non
random possibilities $\cI^\prime\subset {\cal{H}}$  for the set
  $\cI$ and apply Proposition \ref{chi2} in the non random case. As the cardinality of ${\cal{H}}$ is less than $p^L$ by the limitations imposed on $\cI$, we get,
 \begin{align*}
 P\left(\frac{1}{\sigma^2}\|P_{V_{\cI}}[W]\|_{2}^2\geq
 z^2\right)&\leq \sum_{\cI^\prime\subset{\cal{H}}}
 P\left(\frac{1}{\sigma^2}\|P_{V_{\cI^\prime}}[W]\|_{2}^2\geq
 z^2\right)\\
 &\leq p^{L}\exp\left(-z^2/8\right)\\
&\leq
\exp\left(-z^2\left[1/8-\frac{L\log{p}}{z^2}\right]\right)\\
&\leq \exp\left(-z^2/16\right)
 \end{align*}
as soon as $z^2\ge 4\left(\sup\#\{\cI^\prime,\;
\cI^\prime\subset{\cal{H}}\}\right)$ and $z^2\geq 16\;L\log{p}$. To
end up the proof, it remains to observe that $\sup\#\{\cI^\prime,\;
\cI^\prime\subset{\cal{H}}\}\le Lt^*$.

\bibliographystyle{apalike}
\bibliography{BiblioGrol1}

\end{document}